\documentclass[11pt, twoside, leqno]{article}

\usepackage[T1]{fontenc}
\usepackage{amssymb,amsmath,amsthm}
\usepackage{enumerate}
\usepackage[dvips]{graphicx}
\usepackage{epic,eepic}

\newtheorem{lem}{Lemma}
\newtheorem{thm}[lem]{Theorem}
\newtheorem{prop}[lem]{Proposition}
\newtheorem{coro}[lem]{Corollary}

\theoremstyle{definition}

\newtheorem{de}[lem]{Definition}
\newtheorem{ex}[lem]{Example}
\newtheorem{pozn}[lem]{Remark}

\newcommand{\Z}{{\mathbb Z}}
\newcommand{\N}{{\mathbb N}}
\newcommand{\Q}{{\mathbb Q}}
\newcommand{\R}{{\mathbb R}}
\newcommand{\A}{\mathcal{A}}
\newcommand{\B}{\mathcal{B}}

\newcommand{\Scal}{\mathcal{S}}

\newcommand{\coloneq}{\mathrel{\mathop:}=}

\begin{document}


\title{Numbers with integer expansion in the numeration system with negative base}

\author{P. Ambro\v z, D. Dombek, Z. Mas\'akov\'a, E. Pelantov\'a\\[1mm]
Department of Mathematics FNSPE\\ Czech Technical University in Prague\\
Trojanova 13, 120 00 Praha 2, Czech Republic}

\date{}

\maketitle


\renewcommand{\thefootnote}{}

\footnote{2010 \emph{Mathematics Subject Classification}: Primary XXXX; Secondary YYYY.}


\renewcommand{\thefootnote}{\arabic{footnote}}

\setcounter{footnote}{0}


\begin{abstract}
In this paper, we study representations of real numbers in the
positional numeration system with negative basis, as introduced by
Ito and Sadahiro. We focus on the set $\Z_{-\beta}$ of numbers
whose representation uses only non-negative powers of $-\beta$,
the so-called $(-\beta)$-integers. We describe the distances
between consecutive elements of $\Z_{-\beta}$. In case that this
set is non-trivial we associate to $\beta$ an infinite word
$\boldsymbol{v}_{-\beta}$ over an (in general infinite) alphabet.
The self-similarity of $\Z_{-\beta}$, i.e., the property $-\beta
\Z_{-\beta}\subset \Z_{-\beta}$, allows us to find a morphism
under which $\boldsymbol{v}_{-\beta}$ is invariant. On the example of two
cubic irrational bases $\beta$  we demonstrate the difference
between Rauzy fractals generated by $(-\beta)$-integers and by
$\beta$-integers.
\end{abstract}


\section{Introduction}

In~\cite{ItoSadahiro}, Ito and Sadahiro have introduced a new
numeration system, using negative base $-\beta<-1$ for expansion
of real numbers. The method for obtaining the
$(-\beta)$-expansion is analogous to the one for positive bases
introduced by R\'enyi~\cite{Renyi}. During the 50 years since the
publication of this paper, the $\beta$-expansions of R\'enyi were
extensively studied. On the other hand, analogous study of
$(-\beta)$-expansions has yet to be performed.

The R\'enyi expansions of real numbers $x\in[0,1)$ with positive
base $\beta>1$ are constructed using the transformation $T_\beta$
of the unit interval. Simple adaptation leads to the greedy
algorithm which allows one to uniquely expand any positive real
number in the $\beta$-numeration system. Essential tool for
determining which digit strings are admissible in
$\beta$-expansions is the R\'enyi expansion of 1. The criterion
for admissibility was given in~\cite{Parry} using lexicographical
ordering of strings. Frougny and Solomyak~\cite{FruSo}, and then
many others have considered questions about arithmetics on
$\beta$-expansions. Another point of view on $\beta$-numeration is
to study combinatorial properties of the set of positive real
numbers with integer $\beta$-expansion. The so-called
$\beta$-integers were considered for example by Burd\'\i k et
al.~\cite{BuFrGaKr} as natural counting system for coordinates of
points in quasicrystal models. The description of distances
between consecutive $\beta$-integers is due to
Thurston~\cite{Thurston}. Consequently, in case that the R\'enyi
expansion of $1$ is eventually periodic, the sequence of
$\beta$-integers can be coded by an infinite word, denoted
$\boldsymbol{u}_\beta$, over a finite alphabet. Fabre shows that this word is
substitution invariant, see~\cite{Fabre}. A generalization of the
$\beta$-transformation $T_\beta$ is studied from various points of
view in \cite{KaSt}.

Ito and Sadahiro in their paper~\cite{ItoSadahiro} show some
fundamental properties of their new numeration system with
negative basis. Most of them are analogous to $\beta$-expansions,
though often more complicated. For example, the characterization
of admissible digit sequences is given using the alternate order
on digit strings. Ito and Sadahiro also provided a criterion for
the $(-\beta)$-shift to be sofic, and determined explicitly the
absolutely continuous invariant measure of the
$(-\beta)$-transformation $T_{-\beta}$ which is used in the
expansion algorithm.  Frougny and Lai~\cite{ChiaraFrougny} have
studied arithmetical aspects of the $(-\beta)$-numeration. In
particular, they have shown that for $\beta$ Pisot, the
$(-\beta)$-shift is a sofic system and addition is realizable by a
finite transducer.

Our aim is to deepen the knowledge about $(-\beta)$-expansions. We
focus, in particular, on the properties of the set $\Z_{-\beta}$
of $(-\beta)$-integers which are defined in natural analogy to the
classical $\beta$-integers. We describe the bases for which this
set is non-trivial, and we show that $\Z_{-\beta}$ does not have
accumulation points. We try to give an insight to admissibility
and alternate ordering of finite strings in order to describe the
distances between consecutive $(-\beta)$-integers. This is
achieved for a large class of bases $-\beta$. In case that
$\Z_{-\beta}$ is non-trivial we associate to $\beta$ an infinite
word $\boldsymbol{v}_{-\beta}$ over an infinite alphabet.  The
self-similarity of $\Z_{-\beta}$, i.e., the property $-\beta
\Z_{-\beta}\subset \Z_{-\beta}$, allows us to find a morphism
under which $\boldsymbol{v}_{-\beta}$ is invariant. If, moreover, the
corresponding $(-\beta)$-shift is sofic, then we associate an
infinite word $\boldsymbol{u}_{-\beta}$ over a finite alphabet which is
invariant under a primitive morphism. The infinite word
$\boldsymbol{u}_{-\beta}$ is simply constructed from $\boldsymbol{v}_{-\beta}$ by a
letter-to-letter projection. Similar question is studied
in~\cite{Steiner} using a different approach.
 For a Pisot number $\beta$, the set of
$\beta$-integers corresponds naturally to a so-called Rauzy
fractal, see~\cite{akiyama,arnouxito}. Analogous fractal can be
constructed using $(-\beta)$-expansions. For two examples of cubic
Pisot numbers $\beta$ we compare the fractal tiles arising from
$(-\beta)$-expansions with the Rauzy fractal given by the
classical R\'enyi $\beta$-expansions.

\section{$\beta$-expansions versus $(-\beta)$-expansions}\label{sec:porovnani}

\subsection{R\'enyi $\beta$-expansions}

Consider a real base $\beta>1$ and the transformation $T_\beta:[0,1)\to[0,1)$
defined by the prescription $T_\beta(x)\coloneq\beta x-\lfloor \beta x\rfloor$.
The representation of a number $x\in[0,1)$ of the form
\[
x=\frac{x_1}{\beta} + \frac{x_2}{\beta^2} + \frac{x_3}{\beta^3} +\cdots\,,
\]
where $x_i=\lfloor \beta T_\beta^{i-1}(x)\rfloor$, is called the
$\beta$-expansion of $x$. The coefficients $x_i$ are called
digits. We write $\mathrm{d}_\beta(x)=x_1x_2x_3\cdots$ or
$x={\scriptstyle \bullet} x_1x_2x_3\cdots$. If the string
$x_1x_2x_3\cdots$ ends in suffix $0^\omega$, i.e.\ infinite repetition of $0$, we omit it. Since $\beta
T_\beta(x)\in[0,\beta)$, the $\beta$-expansion of $x$ is an
infinite word in the alphabet $\{0,1,\dots,\lceil\beta\rceil-1\}$.

The $\beta$-expansion of an arbitrary real number $x\geq1$ can be
naturally defined in the following way: Find an exponent $k\in\N$
such that $\frac{x}{\beta^k}\in[0,1)$. Using the transformation
$T_\beta$ derive the $\beta$-expansion of $\frac{x}{\beta^k}$ of
the form
\begin{align*}
  \frac{x}{\beta^k} &= \frac{x_1}{\beta} + \frac{x_2}{\beta^2} +
  \frac{x_3}{\beta^3} +\cdots\,,\\
\intertext{so that}
  x &= x_1\beta^{k-1}+x_2\beta^{k-2}+\cdots +
  x_{k-1}\beta+x_k+\frac{x_{k+1}}{\beta}+\cdots\,.
\end{align*}
The $\beta$-expansion of $x$ does not depend on the choice of the
exponent $k$ for which $\frac{x}{\beta^k}\in[0,1)$, which is -- of
course -- not given uniquely. We write $x=x_1x_2\cdots
x_k{\scriptstyle\bullet} x_{k+1}x_{k+2}\cdots$.

The digit string $x_1x_2x_3\cdots$ is said to be
$\beta$-admissible if there exists a number $x\in[0,1)$ so that
$x={\scriptstyle\bullet}x_1x_2x_3\cdots$ is its $\beta$-expansion.
The set of admissible digit strings can be described using the
R\'enyi expansion of 1, denoted by
$\mathrm{d}_\beta(1)=t_1t_2t_3\cdots$, where
$t_1=\lfloor\beta\rfloor$ and
$\mathrm{d}_\beta(\beta-\lfloor\beta\rfloor)=t_2t_3t_4\cdots$. The
R\'enyi expansion of 1 may or may not be finite (i.e., ending in
infinitely many 0's which are often omitted). The infinite R\'enyi
expansion of 1, denoted by $\mathrm{d}_{\beta}^*(1)$ is defined by
\[
\mathrm{d}_{\beta}^*(1) = \lim_{\varepsilon\to0+} \mathrm{d}_\beta(1-\varepsilon)\,,
\]
where the limit is taken over the usual product topology on
$\{0,1,\dots,\lceil\beta\rceil-1\}^\N$. It can be shown that
\[
\mathrm{d}_{\beta}^*(1) =
\begin{cases}
  \bigl(t_1\cdots t_{m-1}(t_m-1)\bigr)^\omega &
  \text{if $\mathrm{d}_\beta(1)=t_1\cdots t_m0^\omega$ with $t_m\neq 0$,} \\
  \mathrm{d}_\beta(1) & \text{otherwise.}
\end{cases}
\]

Parry~\cite{Parry} has shown that the digit string
$x_1x_2x_3\cdots\in\{0,1,\dots,\lceil\beta\rceil-1\}^\N$ is
$\beta$-admissible if and only if for all $i=1,2,3,\dots$
\begin{equation}\label{e:Parryadmissibility}
  0^\omega\preceq_{\text{lex}} x_ix_{i+1}x_{i+2}\cdots
  \prec_{\text{lex}} \mathrm{d}_{\beta}^*(1)\,,
\end{equation}
where $\preceq_{\text{lex}}$ is the lexicographical order.
The lexicographic order of admissible digit strings corresponds to
the standard order of real numbers in $[0,1)$, i.e., if
$\mathrm{d}_\beta(x)=x_1x_2x_3\cdots$ and $\mathrm{d}_\beta(y)=y_1y_2y_3\cdots$ are
$\beta$-expansions of $x$ and $y$ respectively, then $x<y$ if and
only if $\mathrm{d}_\beta(x)\prec_{\text{lex}} \mathrm{d}_\beta(y)$.

Using $\beta$-admissible digit strings, one can define the set of
non-negative $\beta$-integers. According to the knowledge of the
authors, the $\beta$-integers were first defined
in~\cite{BuFrGaKr}. We have
\[
\Z_\beta^+\coloneq\{a_k\beta^k+\cdots + a_1\beta+a_0\mid a_k\cdots
a_1a_00^\omega \text{ is a $\beta$-admissible digit string}\}\,.
\]
In other words, a non-negative real number is a $\beta$-integer,
if its $\beta$-expansion is of the form $x=\sum_{i=0}^k
a_i\beta^i$, i.e., it has no non-zero digits right from the
fractional point ${\scriptstyle\bullet}$.

The distances between consecutive $\beta$-integers are described
in~\cite{Thurston}. It is shown that they take values in the set
$\{\Delta_i \mid i=0,1,2\cdots\}$, where
\begin{equation}\label{eq:dist_thurston}
  \Delta_i=\sum_{j=1}^\infty \frac{t_{i+j}}{\beta^j}\,.
\end{equation}
Since $\Delta_0=1=\frac{t_1}{\beta}+\frac{t_2}{\beta^2}+\frac{t_3}{\beta^3}+\cdots$
and the digit string $t_2t_3t_4\cdots$ satisfies the condition~\eqref{e:Parryadmissibility},
we have $\Delta_i\leq 1$ for all $i=0,1,2,\ldots$.

\subsection{Ito-Sadahiro $(-\beta)$-expansions}\label{sec:ItoSadahiro}

Consider now the real base $-\beta<-1$ and the transformation
$T_{-\beta}: \big[\frac{-\beta}{\beta+1},\frac{1}{\beta+1}\big)
\to \big[\frac{-\beta}{\beta+1},\frac{1}{\beta+1}\big)$
defined by the prescription
$$
T_{-\beta}(x) = -\beta x -\Big\lfloor -\beta x +
\frac{\beta}{\beta+1}\Big\rfloor\,.
$$
Every number
$x\in\big[\frac{-\beta}{\beta+1},\frac{1}{\beta+1}\big)$ can be
represented in the form
\begin{equation}\label{e:rozvojIS}
  x=\frac{x_1}{-\beta} + \frac{x_2}{(-\beta)^2} +
  \frac{x_3}{(-\beta)^3} + \cdots\,,\qquad\text{where}\quad
  x_i=\Big\lfloor -\beta T^{i-1}_{-\beta}(x) + \frac{\beta}{\beta+1}
  \Big\rfloor\,.
\end{equation}
The representation of $x$ in the form~\eqref{e:rozvojIS} is called
the $(-\beta)$-expansion of $x$ and denoted
$$
\mathrm{d}_{-\beta}(x)=x_1x_2x_3\cdots\qquad\text{ or }\qquad
x=0{\scriptstyle\bullet}x_1x_2x_3\cdots\,.
$$
It is shown easily that the digits $x_i$ belong to the set
$\{0,1,\dots,\lfloor\beta\rfloor\}=:\A_\beta$, and thus the string
$x_1x_2x_3\cdots$ belongs to $\A_\beta^\N$. Ito and Sadahiro have
shown that order of reals in $\big[\frac{-\beta}{\beta+1},\frac{1}{\beta+1}\big)$
corresponds to the alternate order $\preceq_{\text{alt}}$ of their $(-\beta)$-expansions,
i.e.,
\[
x<y \quad\iff\quad
\mathrm{d}_{-\beta}(x)=x_1x_2x_3\cdots \prec_{\text{alt}}
\mathrm{d}_{-\beta}(y)=y_1y_2y_3\cdots\,.
\]
Let us recall that the alternate order is defined as follows: We
say that $x_1x_2x_3\cdots \prec_{\text{alt}} y_1y_2y_3\cdots$, if
$(-1)^i(y_i-x_i)>0$ for the smallest index $i$ satisfying $x_i\neq
y_i$.

In order to describe strings that arise as $(-\beta)$-expansions
of some $x\in\big[\frac{-\beta}{\beta+1},\frac{1}{\beta+1}\big)$,
the so-called $(-\beta)$-admissible digit strings, we will use the
notation introduced in~\cite{ItoSadahiro}. Sometimes, we
abbreviate the term $(-\beta)$-admissibility by only admissibility
when no confusion can occur.

We denote $l_{\beta}=\frac{-\beta}{\beta+1}$ and
$r_{\beta}=\frac{1}{\beta+1}$ the left and right end-point of the
domain $I_\beta$ of the transformation $T_{-\beta}$, respectively.
That is $I_\beta=[l_\beta,r_\beta)$. We also denote
\[
\mathrm{d}_{-\beta}(l_\beta)=d_1d_2d_3\cdots\,.
\]

\begin{thm}[\cite{ItoSadahiro}]\label{t:popisretezcu}
  The string $x_1x_2x_3\cdots$ over the alphabet
  $\{0,1,\dots,\lfloor \beta\rfloor\}$ is $(-\beta)$-admissible, if
  and only if for all $i=1,2,3,\dots$,
  \[
  \mathrm{d}_{-\beta}(l_\beta) \preceq_{\text{alt}} x_ix_{i+1}x_{i+2}
  \prec_{\text{alt}} \mathrm{d}_{-\beta}^*(r_\beta)\,,
  \]
  where $\mathrm{d}_{-\beta}^*(r_\beta) = \lim\limits_{\varepsilon\to0+}
  \mathrm{d}_{-\beta}(r_\beta-\varepsilon)$.
\end{thm}

\noindent
The relation between $\mathrm{d}_{-\beta}^*(r_\beta)$ and $\mathrm{d}_{-\beta}(l_\beta)$
is described in the same paper.

\begin{thm}[\cite{ItoSadahiro}]\label{lichy}
  Let $\mathrm{d}_{-\beta}(l_\beta) = d_1d_2d_3\cdots$. If $\mathrm{d}_{-\beta}(l_\beta)$
  is purely periodic with odd period-length, i.e.,
  $\mathrm{d}_{-\beta}(l_\beta) = (d_1d_2\cdots d_{2l+1})^\omega$, then
  $\mathrm{d}_{-\beta}^*(r_\beta) = (0d_1d_2\cdots d_{2l}(d_{2l+1}-1))^\omega$.
  Otherwise, $\mathrm{d}_{-\beta}^*(r_\beta)=0\mathrm{d}_{-\beta}(l_\beta)$.
\end{thm}

\section{$(-\beta)$-expansion of real numbers}\label{sec:minusbetaexpansions}

Analogically to the case of R\'enyi $\beta$-expansions, we use for obtaining
the $(-\beta)$-expansion of an $x\in\R$ a suitable exponent
$l\in\N$ such that
$\frac{x}{(-\beta)^l}\in\big[\frac{-\beta}{\beta+1},\frac{1}{\beta+1}\big)$.
The advantage of $(-\beta)$-expansions is that one can represent
both positive and negative real numbers with the same set of
digits, without using the minus sign, whereas in the case of
R\'enyi expansions, one can represent only non-negative reals. The
disadvantage of $(-\beta)$-expansions is that the choice of the
above mentioned exponent $l\in\N$ may influence the representation
of $x$ as a $(-\beta)$-expansion.

\begin{pozn}
  Consider $x=\frac{\beta^2}{\beta+1}\notin I_\beta$. Since $\frac{x}{-\beta}=l_\beta$,
  we have $\mathrm{d}_{-\beta}(\frac{x}{-\beta})=d_1d_2d_3\cdots$. On the other hand,
  we also have $\frac{x}{(-\beta)^3}=\frac1{-\beta(\beta+1)}\in I_\beta$. In order to find
  its expansion $\mathrm{d}_{-\beta}(\frac{x}{(-\beta)^3})=x_1x_2x_3\cdots$ we compute the
  first digit
  \[
  x_1 = \Big\lfloor -\beta\frac{x}{(-\beta)^3} + \frac{\beta}{\beta+1}\Big\rfloor =
  \Big\lfloor\frac{1}{\beta+1}+\frac{\beta}{\beta+1}\Big\rfloor=1
  \]
  and we have
  $T_{-\beta}\big(\frac{x}{(-\beta)^3}\big)=-\beta\frac{x}{(-\beta)^3}-1=\frac1{\beta+1}-1=l_\beta$.
  Therefore $\mathrm{d}_{-\beta}(\frac{x}{(-\beta)^3})=1d_1d_2d_3\cdots$. This however means that
  \[
  \frac1{-\beta}+\frac{d_1}{(-\beta)^2}+\frac{d_2}{(-\beta)^3}+\frac{d_3}{(-\beta)^4}
  +\cdots = \frac{d_1}{(-\beta)^3} + \frac{d_2}{(-\beta)^4} +
  \frac{d_3}{(-\beta)^5} + \cdots\,.
  \]
\end{pozn}

\begin{lem}\label{l:intervaly}
  Let $a_1a_{2}a_{3}\cdots$ be a $(-\beta)$-admissible digit string with $a_1\neq 0$.
  For fixed $k\in\Z$, denote
  \[
  z=\sum_{i=1}^\infty a_i(-\beta)^{k-i}\,.
  \]
  Then
  \[
  z\in
  \begin{cases}
    \Big[\frac{\beta^{k-1}}{\beta+1},\frac{\beta^{k+1}}{\beta+1}\Big]
     & \text{for $k$ odd},\\[2mm]
    \Big[-\frac{\beta^{k+1}}{\beta+1},-\frac{\beta^{k-1}}{\beta+1}\Big]
     & \text{for $k$ even}.
  \end{cases}
  \]
\end{lem}

\begin{proof}
  Since the alternate order on admissible strings corresponds
  to the order of real numbers, we have
  \begin{equation}\label{e:ner1}
    -\frac{\beta}{\beta+1} \leq \frac{a_1}{-\beta}+
    \frac{a_{2}}{(-\beta)^2}+ \frac{a_{3}}{(-\beta)^3}+ \cdots <
    \frac1{\beta+1}
  \end{equation}
  and
  \begin{equation}\label{e:ner2}
    -\frac{\beta}{\beta+1} \leq \frac{a_{2}}{-\beta}+ \frac{a_{3}}{(-\beta)^2} +
    \frac{a_{4}}{(-\beta)^3}+ \cdots < \frac1{\beta+1}\,.
  \end{equation}

  We find the estimate of the value of $z$. Consider $k$ odd. Then
  multiplying the left inequality of~\eqref{e:ner1} by
  $(-\beta)^{k}=-\beta^k$ we have
  \[
  \frac{\beta^{k+1}}{\beta+1}\geq a_1(-\beta)^{k-1} +
  a_{2}(-\beta)^{k-2}+\cdots =z \,.
  \]
  Multiplying the left inequality of~\eqref{e:ner2} by $(-\beta)^{k-1}=\beta^{k-1}$
  we obtain
  \[
  -\frac{\beta^{k}}{\beta+1}\leq
  a_{2}(-\beta)^{k-2} + a_{3}(-\beta)^{k-3}+\cdots =
  z-a_1(-\beta)^{k-1} =z-a_1\beta^{k-1}\,,
  \]
  which using $a_1\geq 1$ implies
  \[
  z\geq a_1\beta^{k-1}-\frac{\beta^{k}}{\beta+1} \geq \beta^{k-1}
  -\frac{\beta^{k}}{\beta+1}  = \frac{\beta^{k-1}}{\beta+1}\,.
  \]
  Analogous arguments can be used to show the statement for $k$
  even.
\end{proof}

\begin{pozn}\label{pozn:nejednoznacnost}
  It follows from Lemma~\ref{l:intervaly} that the only $z\in\R$ which can
  be represented by two different $(-\beta)$-admissible digit strings are
  \begin{equation}\label{eq:z}
    z=\frac{(-\beta)^k}{\beta+1}\,,\qquad \text{for } k\in\Z\,.
  \end{equation}
  We know that such $z$ belongs to the interval $[l_\beta,r_\beta)$ if and only if $k=1$ or $k\leq{-1}$.
  We have
  \[
  d_{-\beta}\big(\tfrac{(-\beta)^k}{\beta+1}\big) =
  \begin{cases}
    d_1d_2d_3\cdots         & \text{for $k=1$,}\\[1mm]
    0^{-k-1}1d_1d_2d_3\cdots & \text{for $k\leq -1$.}
  \end{cases}
  \]
  Let us stress that even if the digit string $d_1d_2d_3\cdots$ is $(-\beta)$-admissible,
  the string $0d_1d_2d_3\cdots$ is not. In order to keep unicity of $(-\beta)$-expansion of
  all real numbers, we will prefer for numbers $z$ of the form~\eqref{eq:z} the representation using
  the string $1d_1d_2d_3\cdots$, since it is natural to require that admissibility of a string $w$
  implies admissibility of $0w$. Such a convention, however, has an inconvenient consequence:
  For bases $-\beta$, satisfying $1=\frac{(-\beta)^k}{\beta +1}$, the $(-\beta)$-expansion of $1$
  is not equal to $1$. For example, if  $\beta = \frac{1+\sqrt{5}}{2}$ is the golden ratio,
  the preferred $(-\beta)$-expansion of $1$ is equal to $110{\scriptstyle\bullet}$.
\end{pozn}


\section{$(-\beta)$-integers}\label{sec:minusbetaintegers}

In order to avoid ambiguity in defining the expansion of $x$, we
shall define the $(-\beta)$-integers using admissible digit
strings. For a finite digit string $w=a_ka_{k-1}\cdots a_1a_0$
over the alphabet $\A_\beta=\{0,1,\dots,\lfloor\beta\rfloor\}$, we
call its evaluation the value
\[
{\gamma}(w)\coloneq a_k(-\beta)^k+\cdots+a_1(-\beta)+a_0\,.
\]
We define the evaluation of the
empty string $\epsilon$ to be $\gamma(\epsilon)=0$.

\begin{de}
  A real number $x$ is called a $(-\beta)$-integer, if there exists
  a finite string $a_k\cdots a_1a_0$ such that
  \[
  x=\gamma(a_k\cdots a_1a_0)
  \]
  and the digit string $a_k\cdots a_1a_00^\omega$ is
  $(-\beta)$-admissible. The set of all $(-\beta)$-integers is
  denoted by $\Z_{-\beta}$, i.e.,
  \[
  \Z_{-\beta} = \{a_k(-\beta)^k+\cdots+a_1(-\beta)+a_0\mid a_k\cdots
  a_1a_00^\omega\text{ is $(-\beta)$-admissible}\}\,.
  \]
 As explained in Remark~\ref{pozn:nejednoznacnost}, any $(-\beta)$-integer $x$ is the evaluation of an $(-\beta)$-admissible
 string with prefix $0$, which is unique up to the number of prefixed zeros. 
%
\end{de}

\begin{pozn}
  Note that since $0\in I_\beta$ and $T_{-\beta}(0)=0$, we have
  $\mathrm{d}_{-\beta}(0)=0^\omega$. Thus $0\in\Z_{-\beta}$ for every base
  $-\beta$. Another trivial property which follows from
  Theorem~\ref{t:popisretezcu} describing admissible digit strings
  is that $-\beta\Z_{-\beta} \subset \Z_{-\beta}$.
\end{pozn}

\begin{figure}[!ht]
\centering
\setlength{\unitlength}{1.05cm}
\begin{picture}(14.35,1)(0.85,0.5)
\drawline(1.35,1)(13.85,1)
\drawline(1.35,0.9)(1.35,1.1) \put(1.35,0.75){\makebox(0,0)[t]{$z_{-7}$}}
\drawline(2.35,0.9)(2.35,1.1) \put(2.35,0.75){\makebox(0,0)[t]{$z_{-6}$}}
\drawline(2.9,0.9)(2.9,1.1) \put(2.9,0.75){\makebox(0,0)[t]{$z_{-5}$}}
\drawline(3.9,0.9)(3.9,1.1) \put(3.9,0.75){\makebox(0,0)[t]{$z_{-4}$}}
\drawline(4.75,0.9)(4.75,1.1) \put(4.75,0.75){\makebox(0,0)[t]{$z_{-3}$}}
\drawline(5.75,0.9)(5.75,1.1) \put(5.75,0.75){\makebox(0,0)[t]{$z_{-2}$}}
\drawline(6.75,0.9)(6.75,1.1) \put(6.75,0.75){\makebox(0,0)[t]{$z_{-1}$}}
\put(6.75,1.15){\makebox(0.85,0)[b]{$\Delta_1$}}
\drawline(7.6,0.9)(7.6,1.1) \put(7.6,0.75){\makebox(0,0)[t]{$0$}}
\put(7.6,1.15){\makebox(1,0)[b]{$\Delta_0$}}
\drawline(8.6,0.9)(8.6,1.1) \put(8.6,0.75){\makebox(0,0)[t]{$z_1$}}
\put(8.6,1.15){\makebox(0.55,0)[b]{$\Delta_2$}}
\drawline(9.15,0.9)(9.15,1.1) \put(9.15,0.75){\makebox(0,0)[t]{$z_2$}}
\drawline(10.15,0.9)(10.15,1.1) \put(10.15,0.75){\makebox(0,0)[t]{$z_3$}}
\drawline(11,0.9)(11,1.1) \put(11,0.75){\makebox(0,0)[t]{$z_4$}}
\drawline(12,0.9)(12,1.1) \put(12,0.75){\makebox(0,0)[t]{$z_5$}}
\drawline(12.85,0.9)(12.85,1.1) \put(12.85,0.75){\makebox(0,0)[t]{$z_6$}}
\drawline(13.85,0.9)(13.85,1.1) \put(13.85,0.75){\makebox(0,0)[t]{$z_7$}}
\end{picture}
\caption{First few $(-\beta)$-integers for the Tribonacci number $\beta$.}
\label{fig:Tr_minus_int}
\end{figure}
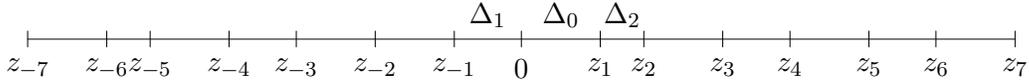

\begin{ex}
  Let $\beta$ be the so-called Tribonacci number, i.e., the real
  root of the polynomial $x^3-x^2-x-1$. Then
  $\mathrm{d}_{-\beta}(l_\beta)=101^\omega$. Let us find the first (in absolute
  value smallest) few $(-\beta)$-integers. Denoting $\Z_{-\beta} =
  \{z_{n}\mid n\in\Z\}$, $z_0=0$, $z_i<z_{i+1}$ for all $i$, we have
  \[
  \begin{aligned}
    1{\scriptstyle\bullet}&=z_1\\
    110{\scriptstyle\bullet}&=z_2\\
    111{\scriptstyle\bullet}&=z_3\\
    100{\scriptstyle\bullet}&=z_4\\
    11011{\scriptstyle\bullet}&=z_5\\
    11000{\scriptstyle\bullet}&=z_6\\
    11001{\scriptstyle\bullet}&=z_7\\
    11110{\scriptstyle\bullet}&=z_8
  \end{aligned}
  \qquad\text{ and }\qquad
  \begin{aligned}
    11{\scriptstyle\bullet}&=z_{-1}\\
    10{\scriptstyle\bullet}&=z_{-2}\\
    1100{\scriptstyle\bullet}&=z_{-3}\\
    1111{\scriptstyle\bullet}&=z_{-4}\\
    1110{\scriptstyle\bullet}&=z_{-5}\\
    1001{\scriptstyle\bullet}&=z_{-6}\\
    1000{\scriptstyle\bullet}&=z_{-7}\\
    1011{\scriptstyle\bullet}&=z_{-8}\,.
  \end{aligned}
  \]
  These points are drawn in Figure~\ref{fig:Tr_minus_int}. The distances between
  consecutive $(-\beta)$-integers take values (cf.~Section~\ref{sec:extremaldigitstrings})
  \[
  z_{i+1}-z_i \in\{\Delta_0=1,\ \Delta_1=\beta-1,\ \Delta_2=\beta^2-\beta-1=\frac{1}{\beta}\}\,.
  \]
\end{ex}

\begin{ex}
  Let $\beta$ be the minimal Pisot number, i.e., the real root of
  the polynomial $x^3-x-1$. In this case
  $\mathrm{d}_{-\beta}(l_\beta)=1001^\omega$. We show that in this case
  $\Z_{-\beta}=\{0\}$. Since $\lfloor\beta\rfloor=1$, the digits in
  admissible strings must belong to the alphabet $\{0,1\}$. If the
  digit string $a_k\cdots a_1a_00^\omega\neq 0^\omega$ was
  $(-\beta)$-admissible, then according to
  Theorem~\ref{t:popisretezcu} also the digit string $10^\omega$
  must be $(-\beta)$-admissible. However, we have $1001^\omega
  \not\preceq_{\text{alt}} 10^\omega$, which is a
  contradiction.
\end{ex}

\noindent The observation from the previous example is generalized
by the following statement, derived from~\cite{MaPeVa}.

\begin{prop}[\cite{MaPeVa}]
  Let $\beta > 1$. Then $\Z_{-\beta}=\{0\}$ if and only if
  $\beta < \frac12 (1 + \sqrt5)$.
\end{prop}

Let us show that the condition of $\beta$ being smaller than the
golden ratio corresponds to the requirement that the
$(-\beta)$-expansion $\mathrm{d}_{-\beta}(l_\beta)$ is of a
special form.

\begin{lem}
  We have $\Z_{-\beta}=\{0\}$ if and only if $10^{2k}1$ is a prefix
  of $\mathrm{d}_{-\beta}(l_\beta)$ for some $k\geq 0$.
\end{lem}
\begin{proof}
  First realize that a $(-\beta)$-admissible string not equal
  to $0^\omega$ exists, if and only if the digit string $10^\omega$
  is $(-\beta)$-admissible. Since $\mathrm{d}_{-\beta}^*(r_\beta)$ starts
  always with 0, the alternate inequality
  $10^\omega\prec_{\text{alt}} \mathrm{d}_{-\beta}^*(r_\beta)$ is always
  satisfied. It can be seen easily, that the other inequality from
  Theorem~\ref{t:popisretezcu}, $\mathrm{d}_{-\beta}(l_\beta)
  \preceq_{\text{alt}} 10^\omega$ is satisfied, if and only if
  no prefix of $\mathrm{d}_{-\beta}(l_\beta)$ has the form $10^{2k}1$ for some
  $k\geq 0$.
\end{proof}

\section{Distances between $(-\beta)$-integers}\label{sec:extremaldigitstrings}

From now on, we suppose that the set $\Z_{-\beta}$ is non-trivial,
i.e., $\beta\geq \frac12(1+\sqrt5)$.

In order to describe distances between adjacent $(-\beta)$-integers, we will study
ordering of finite digit strings in the alternate order. Denote by $\Scal(k)$
the set of infinite $(-\beta)$-admissible digit strings such that erasing a prefix
of length $k$ yields $0^\omega$, i.e., for $k\geq 0$, we have
\[
\Scal(k)=\{a_{k-1}a_{k-2}\cdots a_00^\omega \mid
a_{k-1}a_{k-2}\cdots a_00^\omega \text{ is
$(-\beta)$-admissible}\}\,,
\]
in particular $\Scal(0) = \{0^\omega\}$. For a fixed $k$, the
set  $\Scal(k)$ is finite. Denote by $\mathrm{Max}(k)$ the string
$a_{k-1}a_{k-2}\cdots a_00^\omega$ which is maximal in $\Scal(k)$
with respect to the alternate order and by $\max(k)$ its
prefix of length $k$, i.e., $\mathrm{Max}(k)  = \max(k)0^\omega$.
Similarly, we define $\mathrm{Min}(k)$ and $\min(k)$. Thus,
\[
\mathrm{Min}(k) \preceq_{\text{alt}} r \preceq_{\text{alt}} \mathrm{Max}(k)\,,
\qquad\text{ for all digit strings $r \in \Scal(k)$.}
\]

\begin{pozn}\label{pozn:sufixy}
  For any $a_{k-1}a_{k-2}\cdots a_00^\omega\in{\mathcal S}(k)$ its suffix
  satisfies $a_{j-1}a_{j-2}\cdots a_00^\omega\in{\mathcal S}(j)$ for $j\leq k$.
\end{pozn}

Let us start with a simple observation about the alternate order.

\begin{lem}\label{l:usporadani}
  Let $r^{(1)}$, $r^{(2)}$ be infinite strings over $\A_\beta$ and let $w$ be a finite word over
  $\A_\beta$ of length $|w|$. Then
  \[
  r^{(1)}\preceq_{\text{alt}}r^{(2)} \implies
  \begin{cases}
    wr^{(1)}\preceq_{\text{alt}}wr^{(2)} & \text{if $|w|$ is even,}\\[2mm]
    wr^{(1)}\succeq_{\text{alt}}wr^{(2)} & \text{if $|w|$ is odd.}
  \end{cases}
  \]
\end{lem}

\begin{lem}\label{l:X}
  Let $wdr^{(1)}$, $wcr^{(2)}$ be $(-\beta)$-admissible strings, where the digits $c,d$ satisfy $d>c$.
  Then there exists a digit string $r^{(3)}$ such that $w(d-1)r^{(3)}$ is $(-\beta)$-admissible and
  in the alternate order lies between $wdr^{(1)}$ and $wcr^{(2)}$.
\end{lem}
\begin{proof}
  If $c=d-1$, it suffices to put $r^{(3)}=r^{(2)}$. If $c\leq d-2$, then
  \[
  dr^{(1)}\prec_{\text{alt}} (d-1)r^{(1)}\prec_{\text{alt}} cr^{(2)}\,.
  \]
  If $|w|$ is even, then according to Lemma~\ref{l:usporadani}, we obtain
  \[
  wdr^{(1)}\prec_{\text{alt}} w(d-1)r^{(1)}\prec_{\text{alt}} wcr^{(2)}\,.
  \]
  For $|w|$ odd, we obtain
  \[
  wdr^{(1)}\succ_{\text{alt}} w(d-1)r^{(1)}\succ_{\text{alt}} wcr^{(2)}\,.
  \]
  The admissibility of strings $wdr^{(1)}$, $wcr^{(2)}$ implies that $w(d-1)r^{(1)}$
  is also $(-\beta)$-admissible. Hence we can put $r^{(3)}=r^{(1)}$.
\end{proof}

\begin{lem}\label{l:Y}
  Let $v$ be a finite word over the alphabet $\A_\beta$ and let $vdr^{(1)}$,
  $v(d-1)r^{(2)}$ be $(-\beta)$-admissible strings in ${\mathcal S}(n)$, where
  $n>|v|$. Then the digit strings $vd\mathrm{Min}(k)$ and $v(d-1)\mathrm{Max}(k)$,
  for $k=n-|v|-1$, are $(-\beta)$-admissible strings and no other string from $\Scal(n)$
  lies in between them (with respect to the alternate order).
\end{lem}
\begin{proof}
  We verify admissibility of $vd\mathrm{Min}(k)$ and $v(d-1)\mathrm{Max}(k)$  by showing
  that for every suffix $w$ of $v$ the following inequalities are satisfied
  \[
  \begin{array}{r@{\;}c@{\;}c@{\;}c@{\;}l}
    \mathrm{d}_{-\beta}(l_\beta) & \preceq_{\text{alt}} & wd \mathrm{Min}(k)    &
    \prec_{\text{alt}} & \mathrm{d}_{-\beta}^*(r_\beta)\,,\\[2mm]
    \mathrm{d}_{-\beta}(l_\beta) & \preceq_{\text{alt}} & w(d-1)\mathrm{Max}(k) &
    \prec_{\text{alt}} & \mathrm{d}_{-\beta}^*(r_\beta)\,.
  \end{array}
  \]
  Since $wdr^{(1)}\neq w(d-1)r^{(2)}$ we have one of the following cases.
  \begin{enumerate}[a)]
  \item
    Let $wdr^{(1)}\prec_{\text{alt}} w(d-1)r^{(2)}$. This implies that the length of words $wd$
    and $w(d-1)$ is odd. As according to Remark~\ref{pozn:sufixy} the strings $r^{(1)}$, $r^{(2)}$
    belong to ${\mathcal S}(k)$, we obtain
    \[
    \begin{array}{r@{\;}c@{\;}l@{\quad}c@{\quad}r@{\;}c@{\;}l}
      r^{(1)}&\succeq_{\text{alt}}& \mathrm{Min}(k) & \implies & wdr^{(1)} &\preceq_{\text{alt}}
      &wd \mathrm{Min}(k)\,,\\[2mm]
      r^{(2)}&\preceq_{\text{alt}}& \mathrm{Max}(k) & \implies & w(d\!-\!1)r^{(2)} &\succeq_{\text{alt}}
      &w(d\!-\!1) \mathrm{Max}(k)\,.
    \end{array}
    \]
    Now it suffices to use the transitivity of ordering and admissibility of $wdr^{(1)}$ and
    $w(d-1)r^{(2)}$ to conclude admissibility of $vd\mathrm{Min}(k)$ and $v(d-1)\mathrm{Max}(k)$,
    \[
    \mathrm{d}_{-\beta}(l_\beta)\preceq_{\text{alt}} wdr^{(1)}\preceq_{\text{alt}}
    wd\mathrm{Min}(k)\prec_{\text{alt}}
    w(d-1)\mathrm{Max}(k)\preceq_{\text{alt}} w(d-1)r^{(2)}\prec_{\text{alt}} \mathrm{d}_{-\beta}^*(r_\beta).
    \]
    The latter inequality, together with Lemma~\ref{l:X} already implies that $vd\mathrm{Min}(k)$ and
    $v(d-1)\mathrm{Max}(k)$ are -- in the alternate order -- adjacent strings in $\Scal(n)$
    with $n=|v|+1+k$.
  \item
    Let $wdr^{(1)}\succ_{\text{alt}} w(d-1)r^{(2)}$. In this case words $wd$ and $w(d-1)$ are of even length.
    We obtain
    \[
    w(d-1)r^{(2)}\preceq_{\text{alt}} w(d-1)\mathrm{Max}(k)\prec_{\text{alt}} wd\mathrm{Min}(k)
    \preceq_{\text{alt}} wdr^{(1)}\,,
    \]
    and the argumentation is similar.\qedhere
  \end{enumerate}
\end{proof}

\noindent
The following statement is a direct consequence of Lemma~\ref{l:Y}, taking into account that any pair of $(-\beta)$-integers
is evaluation of a pair of strings of the same length, both starting with 0.

\begin{prop}\label{prop:candidates_for_distances} Let $x <y$ be two consecutive
$(-\beta)$-integers. Then there exist  a unique non-negative integer
$k\in\{0,1,2,\dots\}$ and a positive digit $d \in
\mathcal{A}_\beta\setminus \{0\}$ such that
 words $w(d-1)\mathrm{Max}(k)$ and $wd\mathrm{Min}(k)$ are
$(-\beta)$-admissible strings and
\[
\begin{array}{lcll}
x = \gamma(w(d-1)\max(k)) &<&  y =\gamma(wd\min(k))\quad &\hbox{for $k$ even},\\
x = \gamma(wd\min(k)) &<&  y = \gamma(w(d-1)\max(k)) \quad &\hbox{for $k$ odd},
\end{array}
\]
where $w$ is a finite string over the alphabet $\mathcal{A}_\beta$ with prefix $0$.
In particular,  the distance $y-x$  between these
$(-\beta)$-integers depends only on $k$ and  equals to
\begin{equation}\label{Delty}
\Delta_k:=\Big|(-\beta)^k + \gamma\big(\min(k)\big) -
  \gamma\big(\max(k)\big)\Big|\,.
  \end{equation}
\end{prop}

From the properties of the transformation $T_{-\beta}$ it follows
that the digits $d_i$ of the expansion $\mathrm{d}_{-\beta}(l_\beta)$
satisfy $0\leq d_i\leq d_1$ for all $i$. Imposing more assumptions
on the digits $d_i$ will allow us to describe explicitly the
maximal and minimal strings in ${\cal S}(k)$, and by that also the
distances between consecutive $(-\beta)$-integers.

Note that explicit description of strings $\mathrm{Max}(k)$,
$\mathrm{Min}(k)$ for any given $\beta$ is possible, but providing
an explicit formula for the general case would require very
tedious discussion.

\begin{lem}\label{lem:minmax_infinite}
  Let $\mathrm{d}_{-\beta}(l_\beta)=d_1d_2d_3\cdots$,  where $0<d_i$ and $d_1>d_{2i}$ for all
  $i=1,2,3,\dots$.  Then $\min(0)=\epsilon$, $\max(0)=\epsilon$ and for every $k\geq 0$ we have
  \[
  \begin{aligned}
    \min(2k) &= d_1d_2d_3\cdots d_{2k-1}d_{2k}\,, \\
    \min(2k+1) &= d_1d_2d_3\cdots d_{2k}(d_{2k+1}-1)\,,
  \end{aligned}
  \qquad\text{and}\qquad
    \max(k+1) = 0\min(k)\,.
  \]
\end{lem}

\begin{proof}
  Assumptions on $d_i$ exclude the case that $\mathrm{d}_{-\beta}(l_\beta)$ is purely
  periodic with odd period length. Therefore by Theorem \ref{lichy} one has
  $\mathrm{d}_{-\beta}^*(r_\beta)= 0\mathrm{d}_{-\beta}(l_\beta)$.

  Let us show that $\min(2k)=d_1d_2\cdots d_{2k}$. First we show that the digit string $d_1\cdots
  d_{2k}0^\omega$ is admissible. For that, we need to verify
  inequalities
  \begin{equation}\label{eq:duk1}
    d_1d_2d_3\cdots \preceq_{\text{alt}} d_id_{i+1}\cdots
    d_{2k}0^\omega \prec_{\text{alt}} 0d_1d_2d_3\cdots
  \end{equation}
  for all $i=1,\dots,2k$. By the assumption, we have
  \begin{equation}\label{eq:duk2}
    0^\omega \succ_{\text{alt}} d_{2k+1}d_{2k+2}\cdots\,.
  \end{equation}
  For $i=2r\leq 2k$ the string $d_{2r}d_{2r+1}\cdots d_{2k}$ is of odd length, and thus
  \begin{equation}\label{eq:duk3}
    d_{2r}d_{2r+1}\cdots d_{2k}0^\omega \prec_{\text{alt}} d_{2r}d_{2r+1}\cdots d_{2k}d_{2k+1}d_{2k+2}\cdots
    \prec_{\text{alt}} 0\mathrm{d}_{-\beta}(l_\beta) \,,
  \end{equation}
  where the left inequality is a consequence of~\eqref{eq:duk2} and the right inequality
  follows from admissibility of $\mathrm{d}_{-\beta}(l_\beta)$.
  On the other hand, since $d_{2r}<d_1$, we have directly
  \[
  \mathrm{d}_{-\beta}(l_\beta) \prec_{\text{alt}} d_{2r}d_{2r+1}\cdots d_{2k}0^\omega\,.
  \]
  This, together with~\eqref{eq:duk3} gives~\eqref{eq:duk1} for $i=2r$.

  For $i=2r+1< 2k$ the string $d_{2r+1}\cdots d_{2k}$ is of even
  length and therefore with the use of~\eqref{eq:duk2} and
  admissibility of $\mathrm{d}_{-\beta}(l_\beta)$ itself  we obtain
  \begin{equation}\label{eq:duk4}
    \mathrm{d}_{-\beta}(l_\beta)  \preceq_{\text{alt}} d_{2r+1}\cdots
    d_{2k}d_{2k+1}d_{2k+2}\cdots \prec_{\text{alt}} d_{2r+1}\cdots
    d_{2k}0^\omega\,.
  \end{equation}
  Since $d_{2r+1}>0$, we can claim that $d_{2r+1}\cdots
  d_{2k}0^\omega \prec_{\text{alt}} 0\mathrm{d}_{-\beta}(l_\beta)$.
  Together with~\eqref{eq:duk4}, this confirms validity
  of~\eqref{eq:duk1} for $i=2r+1$.

  It remains to show the minimality of the digit string $d_1\cdots
  d_{2k}0^\omega$ in the set ${\cal S}(2k)$. If there exists a
  string $s$ of length $2k$ such that $d_1\cdots
  d_{2k}d_{2k+1}\cdots \preceq_{\text{alt}} s0^\omega
  \prec_{\text{alt}} d_1\cdots d_{2k}0^\omega$, then, from the
  alternate order, we derive that $d_1\cdots d_{2k}=s$. Thus
  $\min(2k)=d_1d_2\cdots d_{2k}$.

  In order to determine $\min(2k+1)$, we first show that
  $d_1d_2\cdots d_{2k+1}0^\omega$ is not admissible. As
  $d_{2k+2}>0$, we can write that $0^\omega \succ_{\text{alt}}
  d_{2k+2}d_{2k+3}\cdots$, and this implies $d_1d_2\cdots
  d_{2k+1}0^\omega \prec_{\text{alt}} \mathrm{d}_{-\beta}(l_\beta)$.
  Let $\mathrm{Min}(2k+1)=s0^\omega$ for some digit string $s$ of
  length $2k+1$. Since $\min(2k)0^\omega$ belongs to the set ${\cal
    S}(2k+1)$, we have
  \[
  d_1\cdots d_{2k}d_{2k+1}d_{2k+2}\cdots \preceq_{\text{alt}}
  s0^\omega \preceq_{\text{alt}} d_1\cdots d_{2k}0^\omega\,,
  \]
  which implies that $s=d_1d_2\cdots d_{2k}x$ for some digit
  $x\in\A_{\beta}$. Moreover, the digit $x$ is maximal possible, so
  that the string $s0^\omega$ be admissible. It is easy to see that
  $x=d_{2k+1}-1$.

  In order to describe maximal strings, realize that the assumption
  $d_{2i}<d_1$ excludes the possibility that $\mathrm{d}_{-\beta}(l_\beta)$
  is purely  periodic with odd period-length. Therefore
  $\mathrm{d}^*_{-\beta}(r_\beta)=0\mathrm{d}_{-\beta}(l_\beta)$. If for
  some digit string $s$ the string $s0^\omega$ is admissible and not
  equal to $\mathrm{d}_{-\beta}(l_\beta)$ then also $0s0^\omega$ is
  admissible. We immediately obtain $\max(k)=0\min(k-1)$ for $k\geq 1$.
\end{proof}

\begin{thm}\label{thm:mezery_infinite}
  Let $\mathrm{d}_{-\beta}(l_\beta)=d_1d_2d_3\cdots$, where $0<d_{i}$ and $d_1>d_{2i}$ for all
  $i=1,2,3,\dots$.
  Then the distances between adjacent $(-\beta)$-integers take values
  \[
  \Delta_0=1\qquad\text{ and }\qquad\Delta_k=\Big|(-1)^k+\sum_{i=1}^\infty
  \frac{d_{k-1+i}-d_{k+i}}{(-\beta)^i}\Big|\,,\quad k=1,2,3,\dots
  \]
  Moreover, all the distances are less than 2.
\end{thm}

\begin{proof}
  For the description of distances $\Delta_k$  according to \eqref{Delty}
  we need to evaluate $\gamma(\max(k))$ and  $\gamma(\min(k))$.
  By Lemma~\ref{lem:minmax_infinite} we have $\max(k) = 0\min(k-1)$
  for every $k\geq 1$,  and thus $\gamma(\max(k)) =
  \gamma(\min(k-1))$. For the calculation of $\gamma(\min(k))$ we discuss the cases
  of even and odd $k$ separately.

  According to Lemma~\ref{lem:minmax_infinite} we have
  \[
  \begin{split}
    \gamma(\min(2k)) &= \sum_{i=1}^{2k}d_i(-\beta)^{2k-i} =
    \sum_{i=1}^{+\infty}d_i(-\beta)^{2k-i}-
    \sum_{i=2k+1}^{+\infty}d_i(-\beta)^{2k-i} =\\
    &= (- \beta)^{2k} l_\beta - \sum_{i=1}^{+\infty}\frac{d_{2k+i}}{(-\beta)^{i}}\,.
  \end{split}
  \]
  Similarly, we obtain
  \[
  \gamma(\min(2k+1)) = \sum_{i=1}^{2k+1}d_i(-\beta)^{2k+1-i} -1 =
  (- \beta)^{2k+1} l_\beta - 1
  - \sum_{i=1}^{+\infty}\frac{d_{2k+1+i}}{(-\beta)^{i}}\,.
  \]
  Therefore, from~\eqref{Delty},
  \[
  \begin{split}
    \Delta_{2k} &= \Big|(-\beta)^{2k}  + \gamma\big(\min(2k)\big) -
    \gamma\big(\min(2k-1)\big)\Big|= \\
    &= \Big|(-\beta)^{2k}  +  (- \beta)^{2k} l_\beta - (- \beta)^{2k-1}l_\beta +1
    - \sum_{i=1}^{+\infty}\frac{d_{2k+i}}{(-\beta)^{i}} +
    \sum_{i=1}^{+\infty}\frac{d_{2k-1+i}}{(-\beta)^{i}}\Big|\,.
  \end{split}
  \]
  Realizing that $(-\beta)^{2k}  +  (- \beta)^{2k} l_\beta - (-
  \beta)^{2k-1}l_\beta=0$, we obtain the formula for $\Delta_{2k}$.
  The same procedure leads to the description of distances with odd indices.

  It remains to show that $\Delta_k<2$. This follows from the fact that
  $\sum_{i=1}^{+\infty}\frac{d_{r+i}}{(-\beta)^{i}}
  = T^r_{-\beta}(l_\beta) \in I_\beta$, and thus the subtraction of two sums
  in the expression for $\Delta_k$ is equal to the difference of two numbers
  in the interval $I_\beta$, which is of length 1.
\end{proof}

Let us mention that the distances $\Delta_k$ for $k\geq 1$ can be
written in the form
\begin{equation}\label{eq:mezeryvyjadreni}
  \Delta_k = \Big|(-1)^k + T^{k-1}_{-\beta}(l_\beta) -
  T^{k}_{-\beta}(l_\beta)\Big|\,.
\end{equation}
Using such expression, one can describe which distances among
$\Delta_k$ coincide, if $\mathrm{d}_{-\beta}(l_\beta)$ is
eventually periodic. In the following corollary, the lengths $m$
of the pre-period and $p$ of the period are considered the
smallest possible.

\begin{coro}\label{c:splyvanimezeryperiodic}
  Let the digits $d_1,d_2,d_3,\dots$ of
  $\mathrm{d}_{-\beta}(l_\beta)$ satisfy the conditions of
  Theorem~\ref{thm:mezery_infinite}. Then
  \begin{itemize}
  \item
    If $\mathrm{d}_{-\beta}(l_\beta) = d_1\cdots d_m(d_{m+1})^\omega$,
    then $\Delta_{m+k}=\Delta_0$ for all $k\geq 1$.
  \item
    If $\mathrm{d}_{-\beta}(l_\beta) = d_1\cdots d_m(d_{m+1}d_{m+2})^\omega$,
    then $\Delta_{m+k}=\Delta_{m+1}$ for all $k\geq 1$.
  \item
    If $\mathrm{d}_{-\beta}(l_\beta) = d_1\cdots d_m(d_{m+1}\cdots d_{m+p})^\omega$
    and $p$ even, $p\geq 4$, then $\Delta_{m+p+k}=\Delta_{m+k}$ for all $k\geq 1$.
  \item
    If $\mathrm{d}_{-\beta}(l_\beta) = d_1\cdots d_m(d_{m+1}\cdots d_{m+p})^\omega$
    and $p$ odd, $p\geq 3$, then $\Delta_{m+2p+k}=\Delta_{m+k}$ for all $k\geq 1$.
    Moreover, $\Delta_{m+p+k}\neq\Delta_{m+k}$, namely $\Delta_{m+p+k}=2-\Delta_{m+k}$.
  \end{itemize}
\end{coro}

Even though the class of numbers $\beta$ fulfilling the assumption
of Theorem~\ref{thm:mezery_infinite} is quite large, some
interesting cases are omitted; in particular $\beta$ cannot have
finite $\mathrm{d}_{-\beta}(l_\beta)$, i.e.\ with finitely many non-zero digits. The following theorem gives
a result with similar assumption which, moreover, allows
$\mathrm{d}_{-\beta}(l_\beta)$ to be finite. We omit its proof
because it follows the same ideas as the proof of
Lemma~\ref{lem:minmax_infinite}.

\begin{lem}\label{lem:minmax_finite}
  Let $\mathrm{d}_{-\beta}(l_\beta)=d_1d_2\cdots d_m0^\omega$,
  where $d_m\neq 0$ and $m\geq 1$. If $0<d_i$ and $d_1>d_{2i}$ for all $i=1,2,3,4,\dots,m$, then
  \[
  \min(k) =
  \begin{cases}
    d_1\cdots d_k & \text{for even $k<m$,}\\
    d_1\cdots d_{k-1}(d_k-1) &\text{for odd $k<m$,}\\
    d_1\cdots d_m 0^{k-m} & \text{for $k\geq m$,}
  \end{cases}
  \]
  and
  \[
  \max(k) =
  \begin{cases}
    0\min(k-1) & \text{for $k\leq m$,}\\
    0d_1\cdots d_{m-1}(d_m-1) &\text{for even $k=m+1$,}\\
    0d_1\cdots d_{m-1}(d_m+1) &\text{for odd $k=m+1$, $d_m<d_1-1$,}\\
    0d_1\cdots d_{m-2}(d_{m-1}-1)0 &\text{for odd $k=m+1$, $d_m=d_1-1$,}\\
    0d_1\cdots d_{m-1}(d_m+1)\min(1) & \text{for even $k=m+2$,}\\
    0d_1\cdots d_m 0^{k-m-2}1 &\text{for odd $k\geq m+2$,}\\
    0d_1\cdots d_m 0^{k-m-3}1 \min(1) &\text{for even $k\geq m+3$.}
  \end{cases}
  \]
\end{lem}

When using Lemma~\ref{lem:minmax_finite} in formula~\eqref{Delty},
we obtain the following statement describing
the distances $\Delta_k$ in the case of finite
$\mathrm{d}_{-\beta}(l_\beta)$.

\begin{thm}\label{thm:mezery_finite}
  Let $\mathrm{d}_{-\beta}(l_\beta)=d_1d_2\cdots d_m0^\omega$, where $d_m\neq 0$ and $m\geq 1$.
  If $0<d_i$ and $d_1>d_{2i}$ for all $i=1,2,3,4,\dots,m$,
  then the distances between adjacent $(-\beta)$-integers take values
  \begin{align*}
    \Delta_0 &= 1\,,\\
    \Delta_k &= \Big|(-1)^k+\sum_{i=1}^\infty
    \frac{d_{k-1+i}-d_{k+i}}{(-\beta)^i}\Big|\,,\quad
    k=1,\ldots,m-1\,,\\
   \Delta_m &=
   \begin{cases}
     1-\frac{d_m}{\beta}&\text{ for $m$ even}\,,\\
     \frac{d_m}{\beta}&\text{ for $m$ odd}\,,
   \end{cases}\\[1mm]
   \Delta_{m+1}&=
   \begin{cases}
     \Delta_0&\text{ for $m$ even, $d_m<d_1-1$}\,,\\
     \Delta_1&\text{ for $m$ even, $d_m=d_1-1$}\,,\\
     \Delta_0&\text{ for $m$ odd}\,,
   \end{cases}\\[1mm]
   \Delta_{k}&=
   \begin{cases}
     \Delta_0&\text{ for odd $k\geq m+2$}\,,\\
     \Delta_1&\text{ for even $k\geq m+2$}\,.
   \end{cases}
\end{align*}
Moreover, all the distances are less than 2.
\end{thm}
\section{Infinite words associated with $(-\beta)$-integers}

For a positive base $\beta$, if the R\'enyi expansion of unity
$\mathrm{d}_\beta(1)$ is eventually periodic, i.e., $\beta$ is a Parry
number, the ordering of distances between consecutive non-negative
$\beta$-integers can be coded by a right-sided infinite word $\boldsymbol{u}_\beta$
over a finite alphabet.
The word $\boldsymbol{u}_\beta$ is a fixed point of a primitive morphism $\varphi_\beta$,
see~\cite{Fabre} where these morphisms are called canonical
substitutions associated with Parry numbers. The set of
non-negative $\beta$-integers is denoted by $\mathbb{Z}^+_\beta $,
and the set of all $\beta$-integers is defined symmetrically as
$\mathbb{Z}_\beta = \mathbb{Z}^+_\beta \cup
(-\mathbb{Z}^+_\beta)$. Note that besides the symmetry, there is no natural reason for defining
$\mathbb{Z}_\beta$ in this way for $\beta\notin\N$.
The distances in the set of non-positive
$\beta$-integers are then coded by a left-sided infinite word which is invariant
under a morphism different from the canonical substitution $\varphi_\beta$,
which is by a mirror image connected with $\varphi_\beta$.

Let us now concentrate again on the set of $(-\beta)$-integers
which is non-trivial for $\beta$ greater than or equal to the
golden ratio. Note that the set $\Z_{-\beta}$ includes both positive an negative numbers
expressed using non-negative digits.
Our aim is to provide, analogously to the case of
positive base, a morphism fixing a word coding the set of
$(-\beta)$-integers. We start by showing that this is possible for
every $\beta\geq \tau$, irrespectively whether the corresponding
$(-\beta)$-shift is sofic or not, when using an infinite alphabet.

If the $(-\beta)$-shift is sofic, one then finds, by a suitable
projection of the infinite alphabet, a bidirectional infinite word $\boldsymbol{u}_{-\beta}$
over a finite
alphabet fixed by a morphism. The advantage of the negative base number system is that
the left and the right side
of the word $\boldsymbol{ u}_{-\beta}$ are invariant under the same morphism. This is not the case
for positive base system.
On the other hand, unlike the case of positive base systems, where only two types of canonical
substitutions arise according to whether we consider a simple or a
non-simple Parry number, for $(-\beta)$-integers, one does not
obtain a unified prescription for the morphism dependently on the
coefficients $d_i$ of ${\mathrm d}_{-\beta}(l_\beta)$. However, we
give examples of classes of numbers $\beta$ together with
prescriptions for the corresponding morphisms.

Let us recall some notions needed in sequel. If ${\mathcal B}$ is
a (finite or infinite) alphabet, then ${\mathcal B}^*$ denotes the
set of finite words over ${\mathcal B}$. Equipped with the
operation of concatenation, ${\mathcal B}^*$ is a free monoid with
the empty word $\epsilon$ as the neutral element. A morphism over
$\B$ is a mapping $\psi:\B^*\to\B^*$ such that
$\psi(vw)=\psi(v)\psi(w)$ for all pairs of words $v,w\in \B^*$.
The action of a morphism can be naturally extended to infinite
words, both one-directional and pointed bidirectional by
\begin{align*}
\psi(w_0w_1w_2\cdots)&=\psi(w_0)\psi(w_1)\psi(w_2)\cdots\,,\\
\psi(\cdots w_{-2}w_{-1}|w_0w_1w_2\cdots)&=\cdots \psi(w_{-2})
\psi(w_{-1})| \psi(w_0) \psi(w_1)
\psi(w_2)\cdots\,,
\end{align*}
respectively.

Let us assign  to the set $\mathbb{Z}_{-\beta}$ a bidirectional
infinite word $\boldsymbol{ {v}}_{-\beta} = (v_n)_{n\in \mathbb{Z}}$  over
the infinite alphabet $\mathbb{N}$. We show that $\boldsymbol{
{v}}_{-\beta}$ is invariant under a morphism constructed as the
second iteration of an antimorphism $\Phi$. By an antimorphism
over an alphabet $\mathcal{B}$, we understand a mapping  $\varphi:
\mathcal{B}^* \to \mathcal{B}^* $ such that $\varphi(wv)=
\varphi(v)\varphi(w)$ for any $w,v\in \mathcal{B}^* $.

At first, we define the word  $\boldsymbol{ {v}}_{-\beta} = (v_n)_{n\in
\mathbb{Z}}$ associated with $(-\beta)$-integers. Let $(z_n)_{n\in
\mathbb{Z}}$ be the strictly increasing sequence satisfying
$$
 z_0 = 0 \quad \hbox{and} \quad \mathbb{Z}_{-\beta} = \{ z_n
\mid n \in \mathbb{Z}\}\,.
$$
According to Proposition~\ref{prop:candidates_for_distances}, for
any $n \in \mathbb{Z} $ there exist a unique integer $k\in \mathbb{N}$,
a word $w$ with prefix 0, and a letter $d$ such that
\begin{equation}\label{distance}
  z_{n+1}- z_n = \Big|\gamma\bigl(w(d-1)\max(k)\bigr)-
  \gamma\bigl(wd\min(k)\bigr) \Big|\,.
\end{equation}
We put $v_n=k$. Thus the $n$-th letter in the word  $\boldsymbol{
{v}}_{-\beta} $ equals to the maximal   exponent $k$ for  which
the coefficients in the $(-\beta)$-expansions of $z_n$ and
$z_{n+1}$ differ.
The letters of the infinite word $\boldsymbol{ {v}}_{-\beta}$ we have just obtained take values
in the infinite alphabet  $\mathbb{N} = \{0,1,2,\ldots\}$.

The self-similarity of the set $\mathbb{Z}_{-\beta}$, namely  the
property $-\beta \mathbb{Z}_{-\beta} \subset \mathbb{Z}_{-\beta}$,
allows us to find a morphism under which $\boldsymbol{ {v}}_{-\beta}$ is
invariant.

\begin{thm}\label{infinit}
  Let  $\boldsymbol{ {v}}_{-\beta}$ be the word  associated with $(-\beta)$-integers.
  There exists an antimorphism $\Phi: \mathbb{N}^*\to \mathbb{N}^*$ such that
  $\Psi = \Phi^2$ is a non-erasing non-identical morphism and
  $\Psi(\boldsymbol{ {v}}_{-\beta}) = \boldsymbol{ {v}}_{-\beta}$.
\end{thm}
\begin{proof}
  Suppose that $x=z_n< y=z_{n+1}$ are two consecutive
  $(-\beta)$-integers. Obviously, $-\beta x, -\beta y\in-\beta\Z_{-\beta}\subset\Z_{-\beta}$
  and $-\beta y < -\beta x$. Let us study which $(-\beta)$-integers
  lie between $-\beta y$ and $-\beta y$. Our considerations are
  illustrated by Figures~\ref{fig:subst1} and~\ref{fig:subst2}.

  Suppose at first that the distance $z_{n+1}-z_n$ between $x$ and
  $y$ is coded by an even number $v_n=k$. Then
  using~\eqref{distance}, we have
  \begin{equation}\label{eq:xy}
    \begin{aligned}
      x&=\gamma\bigl(w(d-1)\max(k)\bigr)\,,\\
      y&=\gamma\bigl(wd\min(k)\bigr)\,,
    \end{aligned}
    \quad \hbox{and}\quad
    \begin{aligned}
      -\beta x &= \gamma\bigl(w(d-1)\max(k)0\bigr)\,,\\
      -\beta y&=\gamma\bigl(wd\min(k)0\bigr)\,.
    \end{aligned}
  \end{equation}
  Let us realize that
  \[
  d\mathrm{Min}(k) \preceq_{\hbox{\tiny alt}} d\mathrm{Min}(k+1)
  \preceq_{\hbox{\tiny alt}} (d-1)\mathrm{Max}(k+1)
  \preceq_{\hbox{\tiny alt}}(d-1)\mathrm{Max}(k)\,.
  \]
  Therefore according to
  Proposition~\ref{prop:candidates_for_distances}, points
  \begin{equation}\label{eq:x'y'}
    y' = \gamma\bigl( wd\min(k+1)\bigr)\quad \hbox{and}\quad x' =
    \gamma\bigl( w(d-1)\max(k+1)\bigr)
  \end{equation}
  are consecutive $(-\beta)$-integers situated between $-\beta y$
  and $-\beta x$, see Figure~\ref{fig:subst1}. Their distance is
  coded by the number $k+1$ and the value of their distance is
  $\Delta_{k+1}$.
  \begin{figure}
    \begin{center}
      \includegraphics{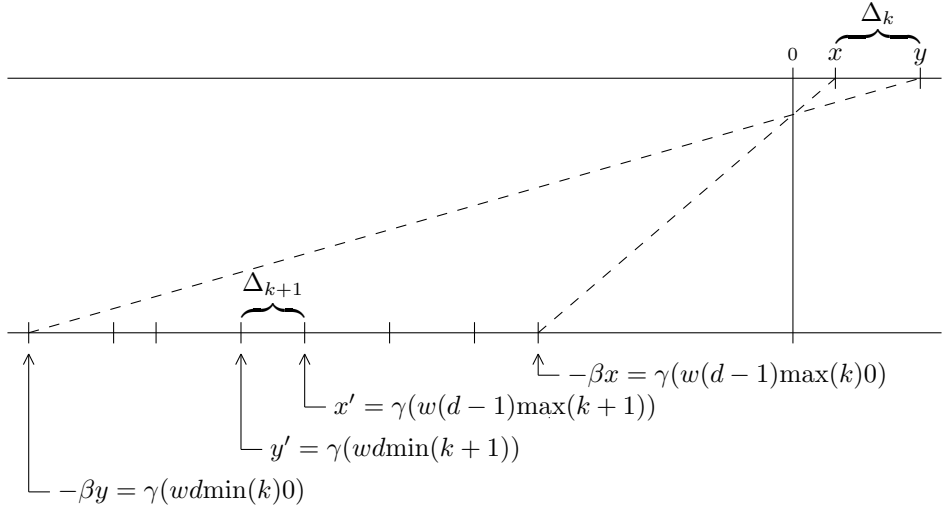}
      \caption{Location of points $x,y$ and $-\beta x,-\beta y$
        from~\eqref{eq:xy}, and $x',y'$ from~\eqref{eq:x'y'}.}
      \label{fig:subst1}
    \end{center}
  \end{figure}

  \noindent For the description of all $(-\beta)$-integers  between
  $-\beta y$ and $-\beta x$, it suffices to determine all digit
  strings $s0^\omega$ such that
  \begin{equation}\label{coMeziMin}
    s0^\omega\in \mathcal{S}(k+1)
    \quad \hbox{and}\quad \min(k+1)0^\omega \preceq_{\hbox{\tiny
        alt}} s0^\omega \preceq_{\hbox{\tiny alt}} \min(k)00^\omega
  \end{equation}
  and to determine all digit strings $r$ such that
  \begin{equation}\label{coMeziMax}
    r0^\omega\in \mathcal{S}(k+1) \quad \hbox{and}\quad
    \max(k)00^\omega \preceq_{\hbox{\tiny alt}} r0^\omega
    \preceq_{\hbox{\tiny alt}} \max(k+1)0^\omega \,.
  \end{equation}

  Note that similar considerations can be carried out in case that
  the distance between $(-\beta)$-integers $x,y$ is coded by an odd
  number $k$. In the following, we consider both $k$ even and odd.

  Let us order the strings $s0^\omega$ satisfying~\eqref{coMeziMin}
  in the alternate order starting with the greatest one and ending
  with the smallest one. We evaluate $\gamma(wds)$ for all of such
  strings $s$. We thus obtain a sequence of consecutive
  $(-\beta)$-integers, say
  \begin{equation}\label{eq:nuch1}
    -\beta y = \gamma\big(wd\min(k)0\big) = b_0,\ b_{1},\ b_{2},\
    \ldots,\ b_{n_k} =\gamma\big(wd\min(k+1)\big) = y'\,.
  \end{equation}
  This sequence is increasing if $k$ is even and decreasing if $k$
  is odd. The distance between $b_{i-1}$ and $b_{i}$, for all
  $i=1,2,\ldots, n_k$, is coded by a number, say $s_i$, in
  $\{0,1,2,\ldots, k\}$. Denote by $S_k$ the word of length $n_k$
  found as the concatenation of $s_1,s_2,\dots, s_{n_k}$,
  \begin{equation}\label{slovoCoMeziMin}
    S_k \coloneq s_1s_2\cdots s_{n_k}\,.
  \end{equation}

  Similarly, we now order the strings $r0^\omega$ satisfying
  \eqref{coMeziMax} starting with the smallest one $\max(k)00^\omega$
  and ending with $\max(k+1)0^\omega$.
  Evaluating $\gamma\big(w(d-1)r\big)$, we obtain a sequence of
  consecutive $(-\beta)$-integers, say
  \begin{equation}\label{eq:nuch2}
    -\beta x = \gamma\big(w(d-1)\max(k)0\big)=c_0,\ c_{1},\ \ldots,\
    c_{m_k} =\gamma\big(w(d-1)\max(k+1)\big) = x'\,.
  \end{equation}
  This sequence is decreasing for $k$ even and increasing for $k$
  odd. Let us denote by $r_i\in \{0,1,2,\ldots, k\}$  the number
  coding the distance between $c_{i-1}$ and $c_{i}$ for
  $i=1,2,\ldots, m_k$ and denote the concatenation
  \begin{equation}\label{slovoCoMeziMax}
    R_k \coloneq r_1r_2\cdots r_{m_k} \,.
  \end{equation}
  Let us stress that words $S_k$ and $R_k$ depend only on $k$ which
  was assigned to the distance between $x=z_{n}$  and $y=z_{n+1}$
  and that $S_k$ and $R_k$ are independent of $w$ and $d$ occurring
  in the evaluation of $z_{n+1}-z_n$ in~\eqref{distance}. The
  situation for even $k$ is depicted at Figure~\ref{fig:subst2}.

  \begin{figure}
    \begin{center}
      \includegraphics{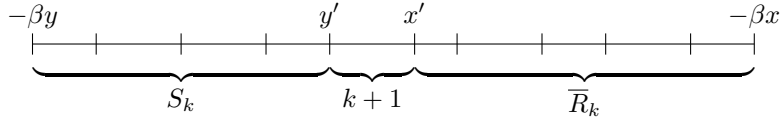}
      \caption{Emplacement of $(-\beta)$-integers $b_0,\dots,b_{n_k}$
        of~\eqref{eq:nuch1} and $c_0,\dots,c_{m_k}$ of~\eqref{eq:nuch2}
        for $k$ even.} \label{fig:subst2}
    \end{center}
  \end{figure}

  Now we are in the position to define the antimorphism $\Phi:
  \mathbb{N}^* \to \mathbb{N}^* $. Put for all $\ell \in \mathbb{N}$
  \begin{equation}\label{anti}
    \Phi(2\ell) = S_{2\ell}\,(2\ell+1)\,\overline{R}_{2\ell}\quad {\rm
      and}\quad \Phi(2\ell +1) =
    R_{2\ell+1}\,(2\ell+2)\,\overline{S}_{2\ell+1}\,.
  \end{equation}
  The notation $\overline{w}$ is used for the mirror image of the
  word $w$, i.e., $\overline{w} = w_tw_{t-1}\ldots w_2w_1$ if
  $w=w_1w_2\ldots w_t$.

  The word $\Phi(k)$, defined in~\eqref{anti}, codes the distances
  between $-\beta y$ and  $-\beta x$ in case that the distance $y-x$
  between consecutive $(-\beta)$-integers $x,y$, $x<y$, is coded by
  $k$.

  \begin{figure}
    \begin{center}
      \includegraphics{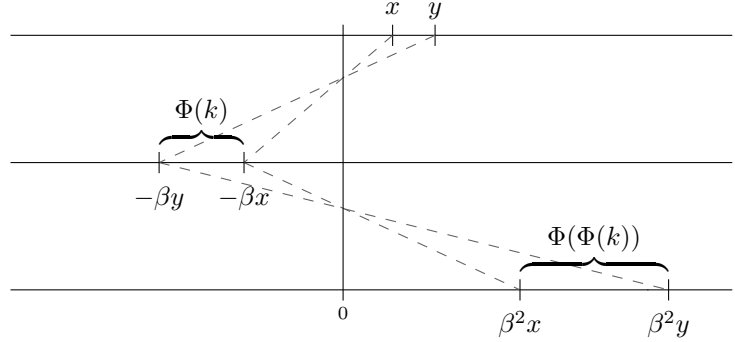}
      \caption{Construction of the morphism $\Psi=\Phi^2$ from the
        antimorphism $\Phi$.} \label{fig:subst3}
    \end{center}
  \end{figure}

  Now applying the same procedure to all pairs of consecutive
  $(-\beta)$-integers occurring between $-\beta y$ and $-\beta x$,
  we can find the word coding the sequence of $(-\beta)$-integers
  between $\beta^2x$ and $\beta^2y$, as shown in
  Figure~\ref{fig:subst3}. Clearly, this word -- depending only on
  $k$ -- is equal to $\Phi^2(k)$.

  To conclude, the self-similarity of $\mathbb{Z}_{-\beta}$
  guarantees that $\boldsymbol{ {v}}_{-\beta}$ is a fixed point of the
  morphism $\Psi=\Phi^2$. The prescription~\eqref{anti} guarantees
  that $\Psi=\Phi^2$ is a non-erasing non-identical morphism.
\end{proof}

Theorem~\ref{infinit} shows that a morphism over an infinite
alphabet, fixing the word $\boldsymbol{v}_{-\beta}$, exists for every
$\beta$ for which $\boldsymbol{v}_{-\beta}$ can be defined. Such morphism
can be explicitly described, whenever strings $\min(k)$ and
$\max(k)$ are known, so that we can determine the words $S_k$
and $R_k$ of~\eqref{anti}. In determining them, we follow the
ideas of the proof of Theorem~\ref{infinit}.

\begin{thm}\label{pronaseD}
  Let the string  $\mathrm{d}_{-\beta}(l_\beta) = d_1d_2d_3\ldots$
  satisfy  $0<d_i$ and $d_1>d_{2i}$ for all $i=1,2,3,\dots$.
  Then the antimorphism from Theorem~\ref{infinit} is of the form
  \begin{align*}
    \Phi(0) &= 0^{d_1-1}1\,,\\
    \Phi(2\ell) &= 0^{d_{2\ell+1}-1}(2\ell+1)0^{d_1-d_{2\ell}-1}1  &&\text{for $\ell\geq 1$,}\\
    \Phi(2\ell+1) &= 0^{d_{2\ell+1}-1}(2\ell+2)0^{d_1-d_{2\ell+2}-1}1 &&\text{for $\ell\geq 0$.}
  \end{align*}
\end{thm}
\begin{proof}
  We need to determine words $S_k$ and $R_k$ for the prescription for
  $\Phi$, given in \eqref{anti}. For that we use the explicit form of strings
  $\min(k)$ and $\max(k)$ from Lemma~\ref{lem:minmax_infinite}.

  By~\eqref{eq:nuch1} and~\eqref{eq:nuch2} we have $S_0=0^{d_1-1}$ and $R_0=\epsilon$,
  thus $\Phi(0) = 0^{d_1-1}1$.
  Next, let us find the word $S_{2\ell}$ with $\ell\geq 1$. Definition
  \eqref{slovoCoMeziMin} of $S_{2\ell}$ requires to find all sequences between
  \[
  \min(2\ell)00^\omega \quad {\rm and} \quad
  \min(2\ell+1)0^\omega\,,
  \]
  or equivalently, to determine all $(-\beta)$-integers between
  \[
  \gamma(d_1d_2 \ldots d_{2\ell}0)\quad {\rm and} \quad
  \gamma(d_1d_2 \ldots d_{2\ell} (d_{2\ell +1}-1))\,.
  \]
  The words \ $d_1d_2 \ldots d_{2\ell}0$  \ and \ $ d_1d_2 \ldots
  d_{2\ell} (d_{2\ell +1}-1)$ \  differ only at the last position
  which in $(-\beta)$-expansion corresponds to the power
  $(-\beta)^0$. Thus all distances between these consecutive
  $(-\beta)$-integers  are coded by $0$,  and therefore  we have
  $S_{2\ell} = 0^{d_{2\ell+1}-1}$. Similarly, $S_0 = 0^{d_1-1}$ as
  $\gamma(\min(0)0) = \gamma(0)$ and $\gamma(\min(1)) =
  \gamma(d_1-1)$.

  Let us determine the words $S_{2\ell+1}$. The complete ordered
  list of all sequences $s0^\omega$, satisfying~\eqref{coMeziMin}
  for $k=2\ell+1$ is
  \[
  \begin{array}{rccccccccl}
    \min(2\ell+1)00^\omega & =
      & d_1 & d_2 & d_3 & \cdots & d_{2\ell} & (d_{2\ell+1}-1) & 0 & 0^\omega \\
    & & d_1 & d_2 & d_3 & \ldots & d_{2\ell} & d_{2\ell+1} & (d_{1}-1)&0^\omega\\
    & & d_1 & d_2 & d_3 & \ldots & d_{2\ell} & d_{2\ell+1} & (d_{1}-2)&0^\omega\\
    & & d_1 & d_2 & d_3 & \ldots & d_{2\ell} & d_{2\ell+1} & (d_{1}-3)&0^\omega\\
    & &     &    &     & \vdots \\
    \min(2\ell+2)0^\omega  & =
      & d_1 & d_2 & d_3 & \ldots & d_{2\ell} & d_{2\ell+1} & d_{2\ell+2}&0^\omega
  \end{array}
  \]
  Therefore we obtain the word $S_{2\ell+1} = 10^{d_1-1-d_{2\ell+2}}$.

 The equality $\max(k)= 0 \min(k-1)$  from Lemma~\ref{lem:minmax_infinite} gives us
immediately $R_{k} = S_{k-1}$ for all $k\geq 1$.  The word $R_0$
is empty, as $\max(0)0^\omega = \max(1)0^\omega=
0^\omega $.
\end{proof}


Let us now study under which conditions one can represent
$(-\beta)$-integers by an infinite word over a restricted finite
alphabet, so that it is still invariant under a primitive
morphism.

\begin{prop}
  Let $\boldsymbol{v}$ be an infinite word over the alphabet $\N$, and let
  $\Psi:\N^*\to\N^*$ be a morphism, such that
  $\Psi(\boldsymbol{v})=\boldsymbol{v}$. Let $\Pi$ be a letter-to-letter morphism
  $\Pi: \mathbb{N}^*\to \mathcal{B}^*$ which satisfies
  \begin{equation}\label{komutaceObecna}
    \Pi \circ\Psi =\Pi \circ \Psi\circ \Pi \,.
  \end{equation}
  Then the infinite word $\boldsymbol{ {u}}=\Pi(\boldsymbol{v})$ is invariant under
  the morphism $\Pi\circ\Psi$.
\end{prop}
\begin{proof}
  We must verify that
  $\big(\Pi\circ\Psi\big)(\boldsymbol{u}) =\boldsymbol{u}$. We write
  \begin{multline*}
    \boldsymbol{ u}=\Pi(\boldsymbol{ {v}}) = \Pi\bigl( \Psi (\boldsymbol{ {v}}))\bigr)
    =\bigl(\Pi\circ \Psi\bigr)(\boldsymbol{ {v}})= \\ =
    \big(\Pi\circ\Psi\circ\Pi\big)(\boldsymbol{ {v}}) = \bigl(\Pi\circ
    \Psi\bigr)\big(\Pi(\boldsymbol{ {v}})\big) =\bigl(\Pi\circ \Psi\bigr)
    (\boldsymbol{ u})\,,
  \end{multline*}
  and hence $\boldsymbol{ u}=\Pi(\boldsymbol{v})$ is a fixed point of the morphism
  $\Pi\circ\Psi$.
\end{proof}

Note that since the morphism $\Psi$ fixing the infinite word
$\boldsymbol{v}_{-\beta}$ coding the set of $(-\beta)$-integers is a power of
an antimorphism $\Phi$, it is sufficient to check that
\begin{equation}\label{komutace}
  \Pi \circ\Phi =\Pi \circ \Phi\circ \Pi \,.
\end{equation}
For we have
\begin{multline*}
  \Pi\bigl( \boldsymbol{ {v}}_{-\beta}\bigr) =
  \Pi\bigl( \Phi (\Phi(\boldsymbol{{v}}_{-\beta} ))\bigr) =
  \Pi\bigl( \Phi (\Pi(\Phi(\boldsymbol{{v}}_{-\beta} )))\bigr) = \\ =
  \Pi\bigl( \Phi (\Pi(\Phi(\Pi(\boldsymbol{{v}}_{-\beta} ))))\bigr) =
  (\Pi\circ \Phi)^2\bigl(\Pi (\boldsymbol{{v}}_{-\beta} )\bigr)\,.
\end{multline*}
Thus the word $\boldsymbol{ {u}}_{-\beta}:=\Pi\big(\boldsymbol{
{v}}_{-\beta}\big)$ is fixed by the morphism
$\psi=\varphi^2=(\Pi\circ \Phi)^2$, where $\varphi= \Pi\circ \Phi$
is an antimorphism over the restricted alphabet $\mathcal{B}$.


Let us now consider the cases of bases where the $(-\beta)$-shift is sofic, i.e.,
such that the $(-\beta)$-expansion $\mathrm{d}_{-\beta}(l_\beta)$ is eventually periodic.
We suggest to call such numbers \emph{Ito-Sadahiro numbers}.

These are the cases where the distances $\Delta_k$ between consecutive
$(-\beta)$-integers take only finitely many values (cf.\ Corollary~\ref{c:splyvanimezeryperiodic}
and Theorem~\ref{thm:mezery_finite}) and thus the set $\Z_{-\beta}$ can be coded by a bidirectional
infinite word $\boldsymbol{ {u}}_{-\beta}$ over a finite alphabet $\mathcal{B}\subset \mathbb{N}$.
Corollary~\ref{c:splyvanimezeryperiodic} and Theorem~\ref{thm:mezery_finite} also suggest
a suitable projection of the infinite alphabet $\N$ to the restricted alphabet ${\mathcal B}$.
By verifying condition~\eqref{komutace}, one can show that $\boldsymbol{ {u}}_{-\beta}$ is
invariant under a primitive morphism. By doing so for eventually periodic
$\mathrm{d}_{-\beta}(l_\beta)=d_1d_2\ldots d_m(d_{m+1}\cdots d_{m+p})^\omega$, one finds that a
prescription in terms of coefficients $d_i$ cannot be written in one formula. Rather it differs
dependently on whether the length of period is shorter or longer, even or odd.
Similarly, it is the case for finite  $\mathrm{d}_{-\beta}(l_\beta)=d_1d_2\ldots d_m0^\omega$.
The discussion is tedious, that is why we show the procedure on only two classes of numbers
$\beta$ together with the corresponding primitive morphism fixing the word
$\boldsymbol{ u}_{-\beta}$.

\begin{ex}\label{ex:1}
  Let $\mathrm{d}_{-\beta}(l_\beta)=d_1d_2\ldots d_m(d_{m+1})^\omega$, where $m\in\N$ is minimal
  possible, i.e., $d_m\neq d_{m+1}$. Assume that for parameters $d_i$ the assumption of
  Theorem~\ref{pronaseD} is satisfied. According to Theorem~\ref{pronaseD}, we have for $m$ even
  \begin{align*}
    \Phi(0) &= 0^{d_1-1}1\,,\\
    \Phi(1) &= 0^{d_{1}-1}20^{d_1-d_{2}-1}1\,,\\
    \Phi(2) &= 0^{d_{3}-1}30^{d_1-d_{2}-1}1\,,\\
    \Phi(3) &= 0^{d_{3}-1}40^{d_1-d_{4}-1}1\,,\\
    &\vdots \\
    \Phi(m) &= 0^{d_{m+1}-1}(m+1)0^{d_1-d_{m}-1}1\,,\\
    \Phi(k) &= 0^{d_{m+1}-1}(k+1)0^{d_1-d_{m+1}-1}1 \qquad\text{for all $k\geq m+1$,}\\
    \intertext{and for $m$ odd}
    \Phi(0) &= 0^{d_1-1}1\,,\\
    \Phi(1) &= 0^{d_{1}-1}20^{d_1-d_{2}-1}1\,,\\
    \Phi(2) &= 0^{d_{3}-1}30^{d_1-d_{2}-1}1\,,\\
    \Phi(3) &= 0^{d_{3}-1}40^{d_1-d_{4}-1}1\,,\\
    &\vdots\\
    \Phi(m) &= 0^{d_{m}-1}(m+1)0^{d_1-d_{m+1}-1}1\,,\\
    \Phi(k) &= 0^{d_{m+1}-1}(k+1)0^{d_1-d_{m+1}-1}1\qquad\text{for all $k\geq m+1$.}
  \end{align*}
  For all $m\in\N$, consider the restricted alphabet $\mathcal{B}=\{0,1,2,\ldots,m\}$
  and the projection $\Pi:\N\to{\mathcal B}$ defined by
  \[
  \Pi(k)=k \quad\text{for $k = 0,1,\ldots, m$}\qquad\text{ and }
  \qquad\Pi(k)=0\quad\text{for $k\geq m+1$.}
  \]
  It can be verified easily that our choice of $\Pi$ guarantees validity of
  condition~\eqref{komutace}.
  The antimorphism $\varphi=\Pi\circ\Phi$ associated with $\beta$ has the form
  \begin{align*}
    \varphi(0) &= 0^{d_1-1}1\,,\\
    \varphi(1) &= 0^{d_{1}-1}20^{d_1-d_{2}-1}1\,,\\
    \varphi(2) &= 0^{d_{3}-1}30^{d_1-d_{2}-1}1\,,\\
    \varphi(3) &= 0^{d_{3}-1}40^{d_1-d_{4}-1}1\,,\\
    &\vdots\\
    \varphi(m) &= 0^{D+d_1-1}1\,,
  \end{align*}
  where $D=(-1)^m(d_{m+1}-d_m)$.

  According to Corollary~\ref{c:splyvanimezeryperiodic}, distances $\Delta_k$ for
  $k=m+1,m+2, m+3, \ldots $ coincide with $\Delta_0=1$. Thus the
  letters projected by $\Pi$ to the letter $0$ code the
  same distance between consecutive $(-\beta)$-integers.
\end{ex}

\begin{ex}\label{ex:2}
  Let $\mathrm{d}_{-\beta}(l_\beta)=d_1d_2\cdots d_m0^\omega$, with $d_m\neq 0$, $m$ even,
  and assume that $0<d_i$ and $d_1>d_{2i}$ for all $i=1,2,3,4,\dots,m$. One finds the
  prescription for the antimorphism $\Phi$ from~\eqref{anti}. Then one considers the
  restricted alphabet
  $\mathcal{B}=\{0,1,2,\ldots,m\}$ and projection $\Pi:\N\to{\mathcal B}$
  \begin{align*}
    \Pi(k) &=k \quad\text{for $k = 0,1,\ldots,m$,} \\
    \Pi(m+1)&=
    \begin{cases}
      0&\text{if }d_m<d_1-1\,,\\
      1&\text{if }d_m=d_1-1\,,
    \end{cases}\\
    \Pi(k)&=\begin{cases}
    0&\text{ for odd $k\geq m+2$,}\\
    1&\text{ for even $k\geq m+2$.}
    \end{cases}
  \end{align*}
  Such projection satisfies~\eqref{komutace} and yields the antimorphism
  \begin{align*}
    \varphi(0) &= 0^{d_1-1}1\,,\\
    \varphi(1) &= 0^{d_{1}-1}20^{d_1-d_{2}-1}1\,,\\
    \varphi(2)&= 0^{d_{3}-1}30^{d_1-d_{2}-1}1\,,\\
    \varphi(3)&= 0^{d_{3}-1}40^{d_1-d_{4}-1}1\,,\\
    &\vdots\\
    \varphi(m) &= 0^{d_1-d_m-1}1.
  \end{align*}
\end{ex}

Similarly as shown in Examples~\ref{ex:1} and~\ref{ex:2}, one can find for every $\beta$
with eventually periodic expansion $\mathrm{d}_{-\beta}(l_\beta)$, a bidirectional infinite
word $\boldsymbol{ u}_{-\beta}$ over a finite alphabet, and a primitive morphism fixing
$\boldsymbol{ u}_{-\beta}$. The word $\boldsymbol{ u}_{-\beta}$ is found by projection of
$\boldsymbol{v}_{-\beta}$ over the alphabet $\N$. Another approach for finding
$\boldsymbol{ u}_{-\beta}$ and the corresponding morphism is used by Steiner
in~\cite{Steiner}. His method uses the notion of the first return map.

\section{$(-\beta)$-integers for two cubic bases and their Rauzy fractals}

In this last section we will demonstrate the results for $\Z_{-\beta}$ on two
particular cubic Pisot numbers $\beta$. The first one is the well-known Tribonacci number
$\beta$, i.e., the real root of $x^3-x^2-x-1$. Note that such $\beta$ does not fulfill
the assumptions of either Theorem~\ref{thm:mezery_infinite} or~\ref{thm:mezery_finite}, and
thus an additional inspection of the distances in $\Z_{-\beta}$ is needed. Nevertheless
this case is very interesting since it bears a striking similarity to the classical
case as demonstrated in the sets $\Z_{-\beta}$ and $\Z_\beta$ as well as on associated
Rauzy fractals.

The second case discussed is the cubic number $\beta>1$, root of
$x^3-2x^2-x+1$. Since $\beta=1+2\cos{\frac{2\pi}{7}}$, it appears
naturally in mathematical models of quasicrystals with 7-fold
symmetry~\cite{BuFrGaKr}. This $\beta$ is an example of a number
covered by Theorem~\ref{thm:mezery_finite}, and, moreover, it
demonstrates that the similarity of $\beta$-  and
$(-\beta)$-numeration observed in the Tribonacci case is not
ubiquitous.

Since we compare the distances between consecutive
$\beta$-integers with distances between consecutive
$(-\beta)$-integers, we differentiate between them by the notation
$\Delta_k^+$ and $\Delta_k^-$, respectively.

\subsection{Tribonacci case}

Let $\beta$ be the Tribonacci number, i.e., the real root
$\beta>1$ of $x^3-x^2-x-1$. Then $\mathrm{d}_\beta(1)=1110^\omega$
and according to~(\ref{eq:dist_thurston}) the distances between
consecutive $\beta$-integers are
\[
\Delta_0^+=1\,,\quad \Delta_1^+=\beta-1\,,\quad\text{and}\quad
\Delta_2^+=\tfrac{1}{\beta}\,.
\]
The infinite word $\boldsymbol{ u}_\beta$ coding the set of
$\beta$-integers is invariant under the canonical substitution
\[
\varphi_\beta:\quad 0\mapsto 01\,,\quad 1\mapsto02\,,\quad
2\mapsto0\,.
\]

Let us inspect the set of $(-\beta)$-integers. According to Section~\ref{sec:ItoSadahiro}
we have
\begin{itemize}
\item
  $\mathrm{d}_{-\beta}(l_\beta)=101^\omega$, admissible digits are $0,1$,
\item
  since $\mathrm{d}_{-\beta}(l_\beta)$ is not purely periodic, we have
  $\mathrm{d}_{-\beta}^*(r_\beta)=0101^\omega$,
\item
  $101^\omega\preceq_{\text{alt}} x_ix_{i+1}x_{i+2}\prec_{\text{alt}} 0101^\omega$ for $i=1,2,3,\dots$
  holds for all $(-\beta)$-admissible digit strings.
\end{itemize}
By Proposition~\ref{prop:candidates_for_distances} one can find the distances in the set
$\Z_{-\beta}$ by evaluating expressions
\[
\left|(-\beta)^l + \gamma\big(\min(l)\big) -
\gamma\big(\max(l)\big)\right| \qquad \text{for $l=0,1,2,\ldots$.}
\]
From the admissibility rule we get the following $\min(l)$
and $\max(l)$: $\min(0)=\epsilon$,
$\min(1)=1$, $\max(0)=\epsilon$ and for every
$k\geq 1$ we have
\[
\begin{aligned}
  \min(2k) &= 10(11)^{k-1}\,, \\
  \min(2k+1) &= 10(11)^{k-1}0\,,
\end{aligned}
\qquad\text{and}\qquad\max(k) = 0\min(k-1)\,.
\]
Having these extremal strings, it can be shown that the distances
between adjacent $(-\beta)$-integers take values
\[
\Delta_0^-=1\,,\quad\Delta_1^-=\beta-1\,,\quad\text{and}\quad\Delta_2^-=\tfrac{1}{\beta}\,.
\]
Therefore in this case the distances between consecutive elements
in $\Z_\beta$ and $\Z_{-\beta}$ coincide. In order to obtain the
morphism under which the infinite word $\boldsymbol{ u}_{-\beta}$ is
invariant, we use Theorem~\ref{infinit}, equation~\eqref{anti}  and the projection
\[
\Pi(k)=k\quad\text{for $k=0,1,2$,\quad and\quad}\Pi(k)=0\quad\text{for $k\geq3$,}
\]
since for such $k$, $\Delta_k^-=\Delta_0^-$. By this, we obtain
the antimorphism
\[
\varphi:\quad 0\mapsto 01\,,\quad 1\mapsto02\,,\quad 2\mapsto0\,.
\]

\begin{figure}[!h]
\begin{minipage}{0.45\textwidth}
  \includegraphics[width=\textwidth]{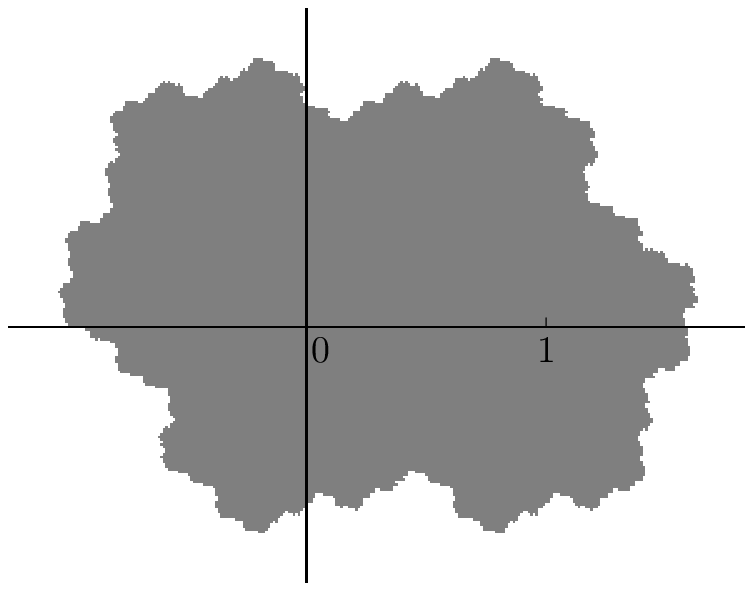}
  \caption{Rauzy fractal for $\Z_\beta$, $\beta$ Tribonacci number.}\label{fig:Tr_Z_beta}
\end{minipage}
\hfill
\begin{minipage}{0.45\textwidth}
  \includegraphics[width=\textwidth]{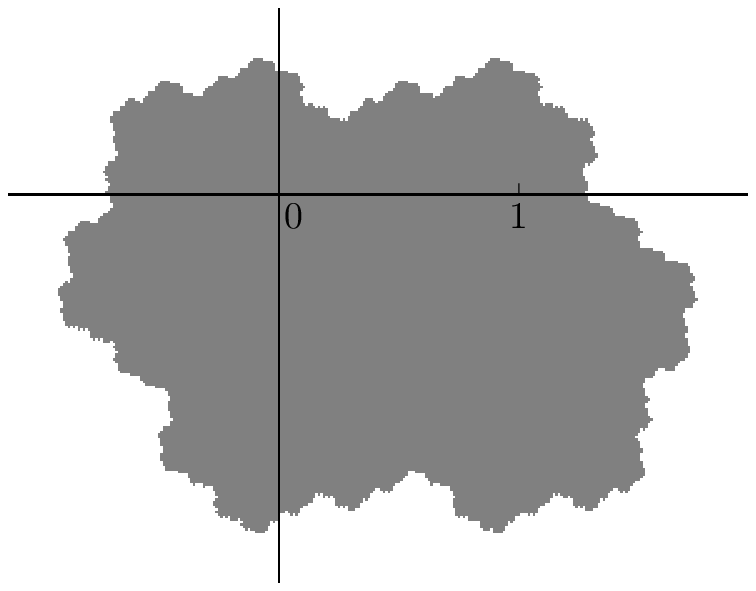}
  \caption{Rauzy fractal for $\Z_{-\beta}$, $\beta$ Tribonacci number.}\label{fig:Tr_Z_-beta}
\end{minipage}
\end{figure}

The similarity between the sets $\Z_\beta$ and $\Z_{-\beta}$ can
be also observed on their Rauzy fractals in
Figures~\ref{fig:Tr_Z_beta} and~\ref{fig:Tr_Z_-beta}. The pictures
represent the set
\[
\Z'_\beta=\{z'\mid z\in\Z_{\beta}\}\,,
\quad\hbox{respectively}\quad \Z'_{-\beta}=\{z'\mid z\in\Z_{-\beta}\}\,,
\]
where $z'$ is the field conjugate of $z\in\Q(\beta)$, i.e.,
\[
\hbox{if }\ z=a+b\beta+c\beta^2\,,\ a,b,c\in\Q\,,\quad\hbox{then}\ z'=a+b{\beta'}+c{\beta'}^2\,,
\]
where $\beta'$ is the conjugate of $\beta$. Since $\beta'$ is a complex number,
$|\beta'|<1$, the sets $\Z'_\beta$, $\Z'_{-\beta}$ are bounded in the complex plane.
The closures of these sets are the so-called Rauzy fractals.

Apparently, the two Rauzy fractals coincide up to a translation in
the complex plane. The reason for this is the following. The
morphism $\varphi^2$ under which $\boldsymbol{ u}_{-\beta}$ is invariant
and the morphism $\varphi_\beta^2$ fixing $\boldsymbol{ u}_\beta$ are
conjugated. Indeed,
\[
\varphi^2: \
\begin{aligned}
  0&\mapsto 0201\\
  1&\mapsto001\\
  2&\mapsto 01
\end{aligned}
\qquad
\varphi_\beta^2: \
\begin{aligned}
  0&\mapsto 0102\\
  1&\mapsto010\\
  2&\mapsto 01
\end{aligned}
\]
and obviously
\[
01\varphi^2(a) = \varphi_\beta^2(a)01 \qquad
\text{for $a\in\{0,1,2\}$.}
\]

\subsection{$\beta$ root of $x^3-2x^2-x+1$}

Let $\beta>1$ be the root of $x^3-2x^2-x+1$. Then
$\mathrm{d}_{\beta}(1)=2(01)^\omega$ and the distances between
consecutive $\beta$-integers are
\[
\Delta_0^+=1\,, \quad \Delta^+_1=\beta-2\,,\quad\text{and}\quad
\Delta_2^+=1-\frac1\beta\,,
\]
and the infinite word $\boldsymbol{ u}_\beta$ is fixed by the canonical
substitution
\[
\varphi_\beta:\quad 0\mapsto001\,,\quad 1\mapsto2\,,\quad
2\mapsto01\,.
\]

The properties of $(-\beta)$-integers are simple to
derive, since $\mathrm{d}_{-\beta}(l_\beta)=210^\omega$ and thus
$\beta$ fulfills the assumptions of
Theorem~\ref{thm:mezery_finite} and Example~\ref{ex:2}.
The distances between adjacent $(-\beta)$-integers are
\[
\Delta_0^- = 1\,,\quad \Delta_1^-=\beta-1\,,\quad\text{and}\quad
\Delta_2^-=1-\frac{1}{\beta}\,.
\]
Clearly, the distances between consecutive elements in $\Z_\beta$
and $\Z_{-\beta}$ are different. Indeed,
$\Delta^-_1=\Delta_1^++1>1$ which cannot happen in the case of
R\'enyi expansions. The infinite word $\boldsymbol{ u}_{-\beta}$ is
invariant under the morphism
$$
\varphi^2:\quad 0\mapsto 02101\,,\quad 1\mapsto 021101\,,\quad
2\mapsto021\,.
$$

\begin{figure}[!h]
\begin{minipage}[b]{0.45\textwidth}
  \includegraphics[width=\textwidth]{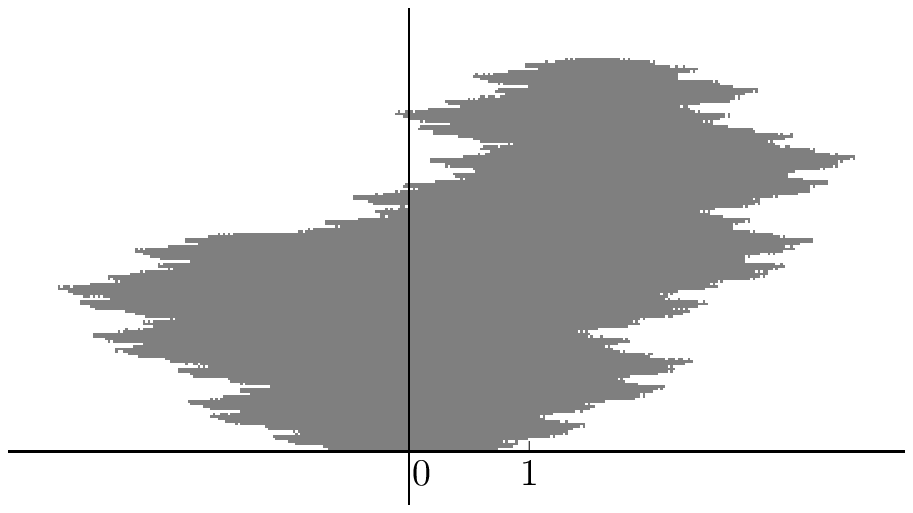}%
  \caption{Rauzy fractal for $\Z_\beta$, $\beta$ root of $x^3=2x^2+x-1$.}%
  \label{fig:Tr_21_beta}
\end{minipage}
\hfill
\begin{minipage}[b]{0.45\textwidth}
  \includegraphics[width=\textwidth]{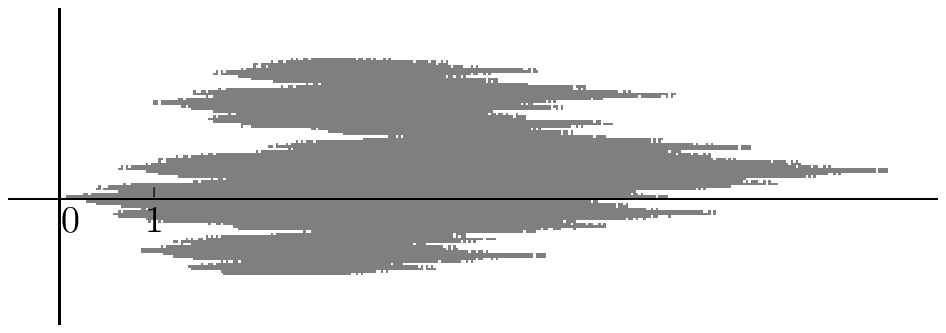}%
  \caption{Rauzy fractal for $\Z_{-\beta}$, $\beta$ root of $x^3=2x^2+x-1$.}%
  \label{fig:Tr_21_-beta}
\end{minipage}
\end{figure}

The dissimilarity of $\Z_\beta$ and $\Z_{-\beta}$ can be also
observed in their Rauzy fractals in Figures~\ref{fig:Tr_21_beta}
and~\ref{fig:Tr_21_-beta}. The base $\beta$ has now two real
conjugates, $\beta'\simeq-0. 8019 $ and $\beta''\simeq 0.5550$. We
consider the sets
$$
\{(x',x'')\mid x\in\Z_\beta\}\,, \quad\hbox{respectively}\quad
\{(x',x'')\mid x\in\Z_{-\beta}\}\,,
$$
where $x',x''$ are images of $x\in\Q(\beta)$ under the isomorphisms
\[
x=a+b\beta+c\beta^2\ \mapsto
\begin{cases}
x'=a+b{\beta'}+c{\beta'}^2\,,\\
x'' = a+b\beta''+c{\beta''}^2\,,
\end{cases}
\]
where $a,b,c\in\Q$.

\section*{Acknowledgements}

We acknowledge financial support by the Czech Science Foundation
grant 201/09/0584 and by the grants MSM6840770039 and LC06002 of
the Ministry of Education, Youth, and Sports of the Czech
Republic. The work was also partially supported by the CTU student
grant SGS10/085/OHK4/1T/14.



\begin{thebibliography}{99}

\bibitem{akiyama}
  S. Akiyama, {\it Pisot number system and its dual tiling},
  In `Physics and Theoretical Computer Science' , ed. by J.P. Gazeau et al.,
  IOS Press (2007), 133--154.

\bibitem{arnouxito}
  P. Arnoux, S. Ito,
  {\it Pisot substitutions and Rauzy fractals},
  Bull. Belg. Math. Soc. Simon Stevin {\bf 8} (2001), 181--207.

\bibitem{BuFrGaKr}
  \v{C}.~Burd\'\i k, Ch.~Frougny, J.~P.~Gazeau, R.~Krejcar,
  {\it Beta-Integers as Natural Counting Systems for Quasicrystals},
  J.~Phys. A: Math. Gen. {\bf 31} (1998), 6449--6472.

\bibitem{Fabre}
  S. Fabre,
  {\it Substitutions et $\beta$-syst\`emes de num\'eration},
  Theoret. Comput. Sci. {\bf 137} (1995), 219--236.

\bibitem{ChiaraFrougny}
  Ch.~Frougny and A.C.~Lai.
  {\it On negative bases},
  In `Proceedings of DLT 09', Lectures Notes in Computer Science {\bf 5583} (2009), 252--263.

\bibitem{FruSo}
  Ch.~Frougny and B.~Solomyak,
  {\it Finite $\beta$-expansions},
  Ergodic Theory Dynam. Systems {\bf 12} (1994), 713--723.


\bibitem{ItoSadahiro}
  S. Ito and T. Sadahiro,
  {\it $(-\beta)$-expansions of real numbers},
  Integers  {\bf 9} (2009), 239--259.


\bibitem{KaSt}
  C.~Kalle, W.~Steiner,
  {\it Beta-expansions, natural extensions and multiple tilings associated with Pisot units},
  to appear in Trans. Amer. Math. Soc., (2010).

\bibitem{MaPeVa}
  Z.~Mas\'akov\'a, E.~Pelantov\'a, T.~V\'avra,
  {\it Arithmetics in number systems with a negative base},
  to appear in Theor. Comp. Sci. (2010), DOI: 10.1016/j.tcs.2010.11.033.

\bibitem{Parry}
  W.~Parry,
  {\it On the $\beta$-expansions of real numbers},
  Acta Math. Acad. Sci. Hung. {\bf 11} (1960), 401--416.

\bibitem{Renyi}
  A. R\'enyi,
  {\it Representations for real numbers and their ergodic properties},
  Acta Math. Acad. Sci. Hung. {\bf 8} (1957), 477--493.

\bibitem{Steiner}
  W.~Steiner, {\it On the structure of $(-\beta)$-integers}, preprint 2010,
  {\tt http://arxiv.org/abs/1011.1755}

\bibitem{Thurston}
  W.P. Thurston,
  {\it Groups, tilings, and finite state automata},
  AMS Colloquium Lecture Notes, American Mathematical Society, Boulder, 1989.

\end{thebibliography}
\end{document}